\documentclass[11pt,a4paper]{article}
\usepackage[margin=2cm]{geometry}
\usepackage[english]{babel}
\usepackage[latin1]{inputenc}
\usepackage{amsmath, amsthm, amssymb}
\usepackage[pdftex]{graphicx}
\usepackage{amssymb}
\usepackage{latexsym}
\usepackage{amstext}
\usepackage{epstopdf}
\usepackage{floatflt}
\usepackage{hyperref}
\usepackage{enumitem}
\usepackage{multirow}
\usepackage{pdflscape}
\usepackage{makecell}
\usepackage{rotating}

\usepackage{eurosym}

\usepackage{slashed}

\usepackage{epstopdf}

\usepackage{nopageno}

\newtheorem{theorem}{Theorem}[section]

\newtheorem{definition}[theorem]{Definition}

\newtheorem{proposition}[theorem]{Proposition}
\newtheorem{remark}[theorem]{Remark}

\newcommand{\adj}{\text{ad}}

\newcommand{\as}{\alpha}

\newcommand{\au}{\mathbf {a_1}}

\newcommand{\ba}{\begin{array}}

\newcommand{\bb}{\beta}

\newcommand{\be}{\begin{equation}}
\newcommand{\bea}{\begin{equation}\begin{array}}
\newcommand{\beal}{\begin{aligned}}
\newcommand{\beas}{\begin{equation*}\begin{array}}
\newcommand{\bef}{\begin{flalign}}
\newcommand{\befs}{\begin{flalign*}}
\newcommand{\bes}{\begin{equation*}}

\newcommand{\bit}{\begin{itemize}}
\newcommand{\blms}{{\mathfrak B}_{{\bf{\mathcal L_{MS}}}}}

\newcommand{\cc}{\mathbb{C}}

\newcommand{\dd}{\mathfrak D}
\newcommand{\ddn}{\mathbf {d_N}}
\newcommand{\ddnd}{\mathbf {d_{N-2}}}
\newcommand{\ddnt}{\mathbf {d_{N-3}}}

\newcommand{\der}{\delta}

\newcommand{\ea}{\end{array}}

\newcommand{\eal}{\end{aligned}}
\newcommand{\eab}{\varepsilon (\alpha,\beta)}
\newcommand{\eag}{\varepsilon (\alpha,\gamma)}
\newcommand{\eba}{\varepsilon (\beta,\alpha)}
\newcommand{\ebg}{\varepsilon (\beta,\gamma)}
\newcommand{\ee}{\end{equation}}
\newcommand{\eea}{\end{array}\end{equation}}
\newcommand{\eeas}{\end{array}\end{equation*}}
\newcommand{\eef}{\end{flalign}}
\newcommand{\eefs}{\end{flalign*}}
\newcommand{\ees}{\end{equation*}}
\newcommand{\ega}{\varepsilon (\gamma,\alpha)}
\newcommand{\egb}{\varepsilon (\gamma,\beta)}

\newcommand{\eit}{\end{itemize}}

\newcommand{\eo}{\mathbf{e_8}}

\newcommand{\eon}{\mathbf{e_8^{(n)}}}

\newcommand{\ep}{\varepsilon}

\newcommand{\es}{\mathbf{e_6}}
\newcommand{\esn}{\mathbf{e_6^{(n)}}}

\newcommand{\est}{\mathbf{e_7}}
\newcommand{\estn}{\mathbf{e_7^{(n)}}}

\newcommand{\fff}{\mathfrak F}

\newcommand{\fq}{\mathbf{f_4}}
\newcommand{\fqn}{\mathbf{f_4^{(n)}}}

\newcommand{\gd}{\mathbf{g_2}}
\newcommand{\gdn}{\mathbf{g_2^{(n)}}}

\newcommand{\gh}{\gamma}

\newcommand{\Ll}{\mathbb L}

\newcommand{\lms}{{\bf{\mathcal L_{MS}}}}

\newcommand{\nin}{\noindent}

\renewcommand{\proof}{\noindent {\bf Proof. }}

\newcommand{\sref}[1]{{\bf\ref{#1}}}

\newcommand{\um}{{\scriptstyle \frac12}}

\newcommand{\zz}{\mathbb Z}

\numberwithin{equation}{section}
\begin{document}

\begin{titlepage}
\begin{center}

\vskip 3.5cm

{\bf \huge Exceptional Periodicity\\\vskip 10pt and Magic Star Algebras\\\vskip 15pt \Large I : Foundations}

\vskip 1.5cm

{\bf \Large Piero Truini${}^{1,2}$, Alessio Marrani${}^{3,4}$, and Michael Rios${}^{5}$ }

\vskip 45pt

{\it ${}^1$ Quantum Gravity Research,\\ 101 S. Topanga Canyon Rd. 1159 Los Angeles, CA 90290 - USA}\\ \vskip 5pt

\vskip 25pt

{\it ${}^2$INFN, sezione di Genova,\\
Via Dodecaneso 33, I-16146 Genova, Italy}

\vskip 25pt

 {\it ${}^3$Museo Storico della Fisica e Centro Studi e Ricerche ``Enrico Fermi'',\\
Via Panisperna 89A, I-00184, Roma, Italy}\\\vskip 5pt

\vskip 25pt

 {\it ${}^4$ Dipartimento di Fisica e Astronomia Galileo Galilei, Universit\`a di Padova,\\and INFN, sezione
di Padova, Via Marzolo 8, I-35131 Padova, Italy}\\\vskip 5pt

\vskip 25pt

{\it ${}^5$Dyonica ICMQG,\\5151 State University Drive, Los Angeles, CA 90032, USA}\\

\vskip 20pt

\texttt{truini@ge.infn.it},
\texttt{alessio.marrani@pd.infn.it},
\texttt{mrios@dyonicatech.com}

\vskip 20pt

\end{center}

\vskip 75pt

\begin{center} {\bf ABSTRACT}\\[3ex]\end{center}


We introduce and start investigating the properties of countably infinite, periodic chains of finite dimensional generalizations of the exceptional Lie algebras: each exceptional Lie algebra (but $\mathbf{g}_{2}$) is part of an infinite family of finite dimensional algebras, which we name \textquotedblleft Magic Star" algebras. These algebras have remarkable similarities with many characterizing features of the exceptional Lie algebras.

\vskip 45pt







%
\vfill

\end{titlepage}

\newpage \setcounter{page}{1} \numberwithin{equation}{section}

\tableofcontents


\section{Introduction}

Operator algebras and their symmetries play a key role in quantum mechanics.
In the attempt to generalize the standard Hilbert space structure of quantum
mechanics over the complex numbers $\cc$,
%
Jordan, Wigner and von Neumann
\cite{JWVN} classified finite
dimensional self-adjoint operator algebras, nowadays named formally real
Jordan algebras. Such a classification singles out the exceptional case of
the Albert algebra (\textit{aka} exceptional Jordan algebra $\mathbf{J}_{%
\mathbf{3}}^{\mathbb{O}}$), related to the octonions $\mathbb{O}$ \cite%
{Octonions} (firstly discovered by J. T. Graves in 1843), which form
the largest normed division algebra\footnote{%
For a review
 of the corresponding Theorem (due to Hurwitz), see \textit{%
e.g.} Th. 1 in \cite{Baez}, and subsequent discussion.}.

On the other hand, Lie algebras have proven to be crucial in the study of
the fundamental interactions of elementary particles, since Gell-Mann's
Eightfold way \cite{GM} and the discovery of the $\mathbf{su}_{3}$ quark and
gluon structures \cite{quark-1, quark-2}. All finite dimensional complex
simple Lie algebras have been classified by Killing and Cartan, and their
non-compact real forms are known: besides the infinite classical series $%
\mathbf{a}_{n}$, $\mathbf{b}_{n}$, $\mathbf{c}_{n}$, $\mathbf{d}_{n}$, five
\textit{exceptional} Lie algebras exist: $\mathbf{g}_{2}$, $\mathbf{f}_{4}$%
, $\mathbf{e}_{6}$, $\mathbf{e}_{7}$ and $\mathbf{e}_{8}$. Along the years,
starting with G\"{u}rsey, Ramond \textit{et al.}, these latter have
characterized the attempts to formulate a Grand Unified Theory of elementary
particles \cite{unif}. In such a framework octonions re-appeared because $%
\mathbf{g}_{2}$, the smallest exceptional Lie algebra, is their algebra of
derivations\footnote{$\mathbf{g}_{2}$ occurs in a number of other physical
contexts, such as, for instance, in the deconfinement phase transitions \cite%
{deco}, in random matrix models \cite{matrix}, in matrix models related to $D
$-brane physics \cite{D-branes}, and in Montecarlo analysis \cite{Montecarlo}%
.}; furthermore, the next largest exceptional Lie algebra, $\mathbf{f}_{4}$,
describes the derivations of the aforementioned Albert algebra \cite{f4}. In
suitable non-compact real forms, all exceptional Lie algebras are used as
electric-magnetic duality ($U$-duality\footnote{%
Here $U$-duality is referred to as the \textquotedblleft
continuous\textquotedblright\ symmetries of \cite{CJ-1}. Their discrete
versions are the $U$-duality non-perturbative string theory symmetries
introduced in \cite{HT-1}.}) algebras in locally supersymmetric theories of
gravity\footnote{%
Some non-compact real forms of exceptional algebras also occur in absence of
local supersymmetry (\textit{cfr. }\cite{Magic-Non-Susy, Romano}, and Refs.
therein).}, and their relation to the Freudenthal-Rozen-Tits Magic Square
was discovered in \cite{MESGT}.

$\mathbf{e}_{7}$ and Lie algebras \textquotedblleft of type $\mathbf{e}_{7}$%
" \cite{brown} have recently appeared in several indirectly related
frameworks of theoretical physics, such as the minimal coupling of vectors
and scalars in cosmology and supergravity \cite{minimal}, in gauge and
global symmetries in the so-called Freudenthal gauge theory \cite{FGT}, and,
by virtue of the so-called black-hole/qubit correspondence (see \cite{BHQIT}
for reviews and list of Refs.), in the entanglement of quantum bits in
quantum information theory.

$\mathbf{e}_{8}$, the largest finite dimensional exceptional Lie algebra,
plays a crucial role in heterotic string theory \cite{Het}, in which the $%
\mathbf{e}_{8}\oplus \mathbf{e}_{8}$ even self-dual lattice corresponds to $%
16$ of the $26$ dimensions of the bosonic string. Moreover, in recent times $%
\mathbf{e}_{8}$ has appeared in other contexts, from mathematics
(computation of the Kazhdan-Lusztig-Vogan polynomials involving $\mathbf{e}%
_{8(8)}$ \cite{KLV}) to experimental physics (namely, in the cobalt niobate
experiment, which is the first actual experiment to detect a phenomenon that
could be modeled using $\mathbf{e}_{8}$ \cite{exp}). After Witten's
formulation of 11-dimensional $M$-theory \cite{9}, the hidden $\mathbf{f}_{4}
$ symmetry of the $D=11$ supermultiplet was observed by Ramond \textit{et al.%
} \cite{10}, and subsequently further studied by Sati \cite{Sati-1, Sati-2}.
Moreover, $\mathbf{J}_{\mathbf{3}}^{\mathbb{O}}$ was also speculated to span
a special charge space related to the 11-dimensional lightcone \cite{11},
since it is naturally endowed with $\mathbf{so}(9)$ and $\mathbf{f}_{4}$
symmetry\footnote{%
Recently, the maximally non-compact (\textit{i.e.}, split) real form $%
\mathbf{f}_{4(4)}$ has been conjectured as the global symmetry of an exotic
ten-dimensional theory in the context of the study of \textquotedblleft
Magic Pyramids" \cite{ICL-Magic, ICL-Magic-2}.}. As far as matrix models are
concerned, the BFSS matrix model \cite{19} for $M$-theory was reformulated
in terms of octonions by Schwarz and Kim in \cite{12}; later on a
Chern-Simons string matrix model was constructed by Smolin exploiting $%
\mathbf{J}_{\mathbf{3}}^{\mathbb{O}}$ \cite{13}, related to Horowitz and
Susskind's conjectured \textquotedblleft bosonic $M$-theory" in $D=27$ \cite%
{14}.

In recent years, advances in algebraic geometry, especially related to works
by Connes and others in the realm of noncommutative geometry \cite{15,16,17}%
, conceived spacetime to be an \textit{emergent} entity, going beyond
Riemannian geometry towards operator algebras. Within this framework, in
which $D$-branes are described by noncommutative coordinates \cite{19,20},
the usual issues with Lorentz symmetry are resolved via discretization,
yielding to an intrinsically fuzzy geometry \cite{18}. Remarkably,
mathematical objects such as $C^{\ast }$-algebras and $K$-homology started
being used in the study of $D$-branes \cite{21, 21-bis}, also determining
the spectral triples of noncommutative geometry and their relevance to the
Standard Model of particle physics \cite{22,23}.

The particular approach to $\mathbf{e}_{8}$ described in a unification model formulated in 2007 by Lisi \cite{L} (later discovered to contain
various issues \cite{DG10}), inspired Truini to rigorously investigate a
special star-like projection - named \textit{\textquotedblleft Magic Star"}
 - of $\mathbf{e}_{8}$ under $\mathbf{a}_{2}$ \cite{Truini}. This led to
a unified construction and characterization of all exceptional Lie algebras,
filling the fourth row of the Freudenthal-Rozen-Tits Magic Square \cite%
{Magic Square}. It was later realized that the Magic Star projection had been
actually envisaged almost ten years before by Mukai, which named it "$%
\mathbf{g}_{2}$ decomposition" \cite{Mukai} and related it to Legendre
varieties. In Truini's formulation, the Magic Star enlightens the structural
relevance of pairs of Jordan algebras of degree three (forming
  \textit{Jordan
pairs}, \cite{loos1}) within each exceptional Lie algebra \cite{Truini}; this
was further investigated in \cite{Marrani-Truini-1, Marrani-Truini-2}, and
also led to some interesting insights in supergravity \cite{Marrani-JP-5}.

Later on, a consistent generalization of exceptional Lie algebras, based on
remarkable properties of the Magic Star projection, \ was introduced by the present
Authors in a contribution to the Proceedings of the 4th Mile High Conference
on Nonassociative Mathematics, held at the University of Denver on July
29-August 5, 2017 \cite{Mile-High-EP}. This resulted in the formulation of
the so-called \textit{"Exceptional Periodicity"}, which generalizes
exceptional Lie algebras to the so-called \textit{\textquotedblleft Magic Star algebras"} \cite%
{Mile-High-EP} parametrized by a
natural number $n\in \mathbb{N}$ (named \textit{\textquotedblleft level"} of Exceptional Periodicity), and enjoying a periodicity (ultimately related to the well known Bott
periodicity). In particular, at each level, the
dimension of the Magic Star algebra is \textit{finite}, raising however the
intriguing question of investigating its $n\rightarrow \infty $ limit; it is
here worth remarking that Exceptional Periodicity was also inspired by the structure of certain $3
$- and $5$- gradings of the exceptional Lie algebras, and especially of $%
\mathbf{e}_{8}$, along with spinor structures\footnote{%
Discussion with Eric Weinstein, during the \textquotedblleft Advances in
Quantum Gravity" symposium, San Francisco, July 2016; see also \cite{Group32-1}}. More details were
presented the year after, in two contributions to the Proceedings of the
32nd International Colloquium on Group Theoretical Methods in Physics, held
in Prague on July 9-13, 2018 \cite{Group32-1, Group32-2}.

The relevance of the Magic Star projection and of Jordan Pairs in the mathematical
description of the fundamental interactions of elementary particles, as
well as for an axiomatic formulation of a consistent theory of quantum
gravity, was started to be investigated in \cite{Truini}, and subsequently
discussed in \cite{Marrani-Truini-Interactions} and in \cite{Group32-2};
recently, a quantum model for the universe at its early stages (including a
mechanism for the creation of space), starting from an initial quantum state
and driven by $\mathbf{e}_{8}$ interactions, was presented by Truini in \cite%
{Truini-EUG}.

\bigskip

The aim of the present paper, which is the first of a series, is to
rigorously establish the mathematical formulation of Exceptional Periodicity\footnote{%
The relevance of Exceptional Periodicity to super Yang-Mills theories in higher dimensions, as
well as to $M$-theory, bosonic string theory and Monstrous\ CFT, has been
recently discussed in \cite{Chester-1, Chester-2}. }. We will prove the
existence of \textit{periodic infinite chains of finite dimensional generalisations of
the exceptional Lie algebras}. In particular, $\mathbf{e}_{8}$ will be shown to be part of an
countably infinite family of algebras (named Magic Star algebras), which resemble
lattice vertex algebras.
Remarkably, for $n=1$, the star-shaped (or $\mathbf{g}_{2}$ decomposition, as
Mukai worded it) structure of known finite dimensional exceptional Lie
algebras is recovered.

As it has already been pointed out in \cite{Mile-High-EP, Group32-1,
Group32-2}, it should be remarked that a key feature of Magic Star algebras is that
they are finite dimensional but \textit{not }of Lie type, namely that they
will not satisfy Jacobi identities anymore; this comes with no surprise,
since Cartan-Killing classification yields that no finite dimensional
exceptional Lie algebras larger than $\mathbf{e}_{8}$ exist. Indeed, within Exceptional Periodicity we will not be dealing with root systems, but rather with \textit{%
\textquotedblleft extended"} root systems, which will be thoroughly defined
further below.

In this perspective, it can thus be stated that Exceptional Periodicity provides a way to go beyond $\mathbf{e}_{8}$ which is radically
different from the way provided by affine and (extended) Kac-Moody
generalizations, such as\footnote{%
For recent development on $\mathbf{e}_{11}$ and beyond, see \cite{beyond e11}%
, and also \cite{Truini-QGR-forthcoming}.} $\mathbf{e}_{8}^{+}=:\mathbf{e}%
_{9}$, $\mathbf{e}_{8}^{++}=:\mathbf{e}_{10}$, $e_{8}^{+++}=:\mathbf{e}_{11}$%
, which also appeared as symmetries for (super)gravity models reduced to $%
D=2,1,0$ dimensions (see e.g. \cite{extended-Refs, West}), respectively, as well as
near spacelike singularities in supergravity \cite{spacelike-Refs}.
\\

Indeed, while such extensions of $\mathbf{e}_{8}$ are still of Lie type but \textit{%
infinite} dimensional \cite{kac}, Magic Star algebras are \textit{not }of Lie type, nevertheless they are \textit{finite} dimensional, for each level of the Exceptional Periodicity itself. Moreover, they maintain the same structure of the finite dimensional exceptional Lie algebras with respect to their maximal orthogonal Lie subalgebra with its spinor representations. We would also like to stress that the product of Magic Star algebras is antisymmetric, and it does not satisfy the Jacobi identity only in the case in which all three entries lie in the spinorial sector. This is the price to pay to have a non-trivial algebraic structure in such a sector. Furthermore, Magic Star algebras provide a generalization of cubic Jordan algebras, consisting in rank-3 matrix algebras introduced by Vinberg in \cite{Vinberg}; we will investigate this interesting issue in a forthcoming paper \cite{EP4}.

\bigskip

The paper is organized as follows.

In Sec. \sref{s:sr} we recall the standard parametrization of the root system of finite dimensional exceptional Lie algebras. Then, in Sec.\sref{s:EPgr} we define the generalized roots of Exceptional Periodicity, by enforcing a Bott periodicization on a type of representation theoretical decomposition of exceptional Lie algebras which highlights their spinorial content (anticipated in \cite{Mile-High-EP,
Group32-2}); in particular, we determine the grouping and content of such generalized roots under suitable two-dimensional projections, which keep the same star-shaped structure of the aforementioned Magic Star projection of exceptional Lie algebras, while generalizing it in a Bott-periodic and infinitely numerable way (which keeps finite dimensionality at each level); such a grouping of the generalized roots is reported in five Tables at the end of the paper. The resulting Magic Star algebras are introduced in Sec. \sref{s:lms} and the properties of the asymmetry function, which is crucial for their definition, are investigated in Sec. \sref{s:epsilon}. Finally, the derivations and automorphisms of the Magic Star algebras are studied in Sec. \sref{s:da}, and proved to be given by their orthogonal Lie algebraic component. Some further developments are recalled and summarized in the concluding considerations in Sec. \sref{s:outlook}.\\


\section{The standard roots}\label{s:sr}

Let $V$ be a Euclidean space of dimension $R$ and $\{k_1 , ... , k_R\}$ an orthonormal basis in $V$. A standard way of writing the roots of the exceptional Lie algebras is the following, \cite{bour}:

\nin
$\gd$ (12 roots)
\bea{ll}
\pm (k_i - k_j) &1\le i<j\le 3\\
\pm \frac13 (-2 k_i + k_{i+1} + k_{i+2}) & i=1,2,3 \text{ (mod 3)}
\eea

\nin
$\fq$ (48 roots)
\bea{ll}
\pm k_i &i=1,...,4\\
\pm k_i\pm k_j &1\le i<j\le 4\\
\frac12 (\pm k_1 \pm k_2 \pm k_3 \pm k_4)
\eea

\nin
$\es$ (72 roots)
\bea{ll}
\pm k_i\pm k_j & 1\le i<j\le 5\\
\frac12 (\pm k_1 \pm k_2 \pm k_3 \pm k_4 \pm k_5 \pm \sqrt{3}k_6) &\text{even \# of +}
\eea

\nin
$\est$ (126 roots)
\bea{ll}
\pm \sqrt{2} k_7\\
\pm k_i\pm k_j & 1\le i<j\le 6\\
\frac12 (\pm k_1 \pm k_2 \pm k_3 \pm k_4 \pm k_5 \pm k_6 \pm \sqrt{2}k_7) &\text{even \# of }+\frac12
\eea

\nin
$\eo$ (240 roots)
\bea{ll}
\pm k_i\pm k_j & 1\le i<j\le 8\\
\frac12 (\pm k_1 \pm k_2 \pm k_3 \pm k_4 \pm k_5 \pm k_6 \pm k_7 \pm k_8) &\text{even \# of +}
\eea

Note that the roots of $\gd$ and $\fq$ can be obtained from those of $\eo$ by respectively projecting on the plane spanned by $k_1-k_2$ and $k_1+k_2-2k_3$, and on the 4-dimensional space spanned by $k_1,k_2,k_3,k_4$.

\nin Moreover one can write the roots of $\es$ and $\est$ as a subset of those of $\eo$ as follows:

\nin
$\es$ (72 roots)
\bea{ll}
\pm k_i\pm k_j & 1\le i<j\le 5\\
\frac12 (\pm k_1 \pm k_2 \pm k_3 \pm k_4 \pm k_5 \pm (k_6+k_7+k_8)) &\text{even \# of +}
\eea

\nin
$\est$ (126 roots)
\bea{ll}
\pm k_7 +k_8\\
\pm k_i\pm k_j & 1\le i<j\le 6\\
\frac12 (\pm k_1 \pm k_2 \pm k_3 \pm k_4 \pm k_5 \pm k_6 \pm (k_7+k_8)) &\text{even \# of +}
\eea


\section{Generalized roots and Exceptional Periodicity}\label{s:EPgr}

We now introduce \textquotedblleft generalized" roots, which do not obey the Weyl reflection symmetry, nor that $2\dfrac{(\alpha,\beta)}{(\alpha,\alpha)}$ be integer for all roots $\alpha$, $\beta$.

For any $n=1,2,...$ we denote $N:=4(n+1)$ and define the generalized roots of $\eon$ as:

\nin
$\eon$ :
\bea{lllcl}
\pm k_i\pm k_j & 1\le i<j\le N && 2N(N-1) &\text{roots}\\
\frac12 (\pm k_1 \pm k_2 \pm ... \pm k_N) &\text{even \# of +} && 2^{N-1} &\text{roots}
\eea

Note that this is a root system only in the case $n=1$, being $\eo^{(1)}=\eo$.

The \textquotedblleft generalized" roots of $\gdn$, $\fqn$, $\esn$, $\estn$ are then obtained by those of $\eon$ in a way similar to the one discussed at the end of Sec. \sref{s:sr}. The {\it generalized} roots of $\gdn$ are obtained by projecting on the space spanned by $k_1-k_2$ and $k_1+k_2-2k_3$, hence $\gdn=\gd$ (i.e., $\gd$ is not generalized in Exceptional Periodicity); those of $\fqn$ are the projection on the space  spanned by $k_1,k_2,...,k_{N-4}$:

\nin
$\fqn$:
\bea{ll}
\pm k_i \ , \ \pm k_i\pm k_j & 1\le i<j\le N-4\\
\frac12 (\pm k_1 \pm k_2 \pm ... \pm k_{N-4}) &
\eea

and finally:\\

\nin
$\esn$:
\bea{ll}
\pm k_i\pm k_j & 1\le i<j\le N-3\\
\frac12 (\pm k_1 \pm k_2 \pm ... \pm (k_{N-2}+k_{N-1}+k_N)) &\text{even \# of +}
\eea

\nin
$\estn$:
\bea{ll}
\pm(k_{N-1}+k_N)\\
\pm k_i\pm k_j & 1\le i<j\le N-2\\
\frac12 (\pm k_1 \pm k_2 \pm ... \pm (k_{N-1}+k_N)) &\text{even \# of +}
\eea

\begin{figure}[htbp]\centering
\includegraphics[width=80mm]{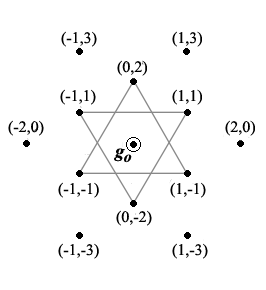}
	\caption{The Magic Star, in the plane coordinatized by $(r,s)$}\label{fig:ms}
\end{figure}

\begin{table}[h]
\begin{equation*}
\begin{array}{| c l c | c | c |}
\hline
\hspace{2.9cm} generalized\ roots & & (r,s) & \#\ of\ roots\\
\hline
\pm (k_1 -  k_2) &  & \pm (2,0) & 2\\
\pm (k_2 -  k_3) &  & \pm (-1,3) & 2\\
\pm (k_3 -  k_1) & & \pm (-1,-3) & 2\\
\hline
&&&\\
 \pm k_i\ , \ \pm k_i \pm  k_j & 4 \le i < j \le N-4 &  & 2(N-7)^2\\
&&(0,0)&\\
\frac{1}{2} (\pm (k_1 + k_2 + k_3) \pm k_4 \pm ... \pm k_{N-4}) & & & 2^{N-6} \\
&&&\\
\hline
&&&\\
k_1 +  k_2 \ , \  - k_3\ , \  - k_3 \pm k_i & i=4,... , N-4 &  & 2N-12\\
&&(0,2)&\\
\frac{1}{2} (k_1 + k_2 - k_3 \pm k_4 \pm ... \pm k_{N-4}) & &  & 2^{N-7} \\
&&&\\
\hline
&&&\\
- k_1 -  k_2 \ , \  k_3\ , \  k_3 \pm k_i & i=4,... , N-4 &  & 2N-12\\
&&(0,-2)&\\
\frac{1}{2} (- k_1 - k_2 + k_3 \pm k_4 \pm ... \pm k_{N-4}) & & & 2^{N-7} \\
&&&\\
\hline
&&&\\
- k_2 -  k_3 \ , \  k_1  \ , \  k_1 \pm k_i & i=4,... , N-4 &  & 2N-12\\
&&(1,1)&\\
\frac{1}{2} (k_1 - k_2 - k_3 \pm k_4 \pm ... \pm k_{N-4}) &  &  & 2^{N-7} \\
&&&\\
\hline
&&&\\
k_2 +  k_3 \ , \ - k_1 \ , \ - k_1 \pm k_i & i=4,... , N-4 &  & 2N-12\\
&&(-1,-1)&\\
\frac{1}{2} (- k_1 + k_2 + k_3 \pm k_4 \pm ... \pm k_{N-4}) & & & 2^{N-7} \\
&&&\\
\hline
&&&\\
- k_1 -  k_3 \ , \  k_2 \ , \  k_2 \pm k_i & i=4,... , N-4 &  & 2N-12\\
&&(-1,1)&\\
\frac{1}{2} (- k_1 + k_2 - k_3 \pm k_4 \pm ... \pm k_{N-4}) & &  & 2^{N-7} \\
&&&\\
\hline
&&&\\
k_1 +  k_3 \ , \  - k_2 \ , \  - k_2 \pm k_i & i=4,... , N-4 &  & 2N-12\\
&&(1,-1)&\\
\frac{1}{2} (k_1 - k_2 + k_3 \pm k_4 \pm ... \pm k_{N-4}) & & & 2^{N-7} \\
&&&\\
\hline
\end{array}
\end{equation*}
\vspace{-16pt}
\caption{The Magic Star for $\fqn$}\label{t:magicfq}
\end{table}

\begin{table}[h]
\begin{equation*}
\begin{array}{| c l c | c | c |}
\hline
\hspace{2.9cm} generalized\ roots & & (r,s) & \#\ of\ roots\\
\hline
\pm (k_1 -  k_2) &  & \pm (2,0) & 2\\
\pm (k_2 -  k_3) &  & \pm (-1,3) & 2\\
\pm (k_3 -  k_1) & & \pm (-1,-3) & 2\\
\hline
&&&\\
 \pm k_i \pm  k_j & 4 \le i < j \le N-3 &  & 2(N-6)(N-7)\\
&&(0,0)&\\
\frac{1}{2} (\pm (k_1 + k_2 + k_3) \pm k_4 \pm ... \pm {\bf u}) & \text{even \# of +} & & 2^{N-5} \\
&&&\\
\hline
&&&\\
k_1 +  k_2 \ , \  - k_3 \pm k_i & i=4,... , N-3 &  & 2N-11\\
&&(0,2)&\\
\frac{1}{2} (k_1 + k_2 - k_3 \pm k_4 \pm ... \pm {\bf u}) & \text{even \# of +} &  & 2^{N-6} \\
&&&\\
\hline
&&&\\
- k_1 -  k_2 \ , \  k_3 \pm k_i & i=4,... , N-3 &  & 2N-11\\
&&(0,-2)&\\
\frac{1}{2} (- k_1 - k_2 + k_3 \pm k_4 \pm ... \pm {\bf u}) & \text{even \# of +} & & 2^{N-6} \\
&&&\\
\hline
&&&\\
- k_2 -  k_3 \ , \  k_1 \pm k_i & i=4,... , N-3 &  & 2N-11\\
&&(1,1)&\\
\frac{1}{2} (k_1 - k_2 - k_3 \pm k_4 \pm ... \pm {\bf u}) & \text{even \# of +} & & 2^{N-6} \\
&&&\\
\hline
&&&\\
k_2 +  k_3 \ , \ - k_1 \pm k_i & i=4,... , N-3 &  & 2N-11\\
&&(-1,-1)&\\
\frac{1}{2} (- k_1 + k_2 + k_3 \pm k_4 \pm ... \pm {\bf u}) & \text{even \# of +} & & 2^{N-6} \\
&&&\\
\hline
&&&\\
- k_1 -  k_3 \ , \  k_2 \pm k_i & i=4,... , N-3 &  & 2N-11\\
&&(-1,1)&\\
\frac{1}{2} (- k_1 + k_2 - k_3 \pm k_4 \pm ... \pm {\bf u}) & \text{even \# of +} & & 2^{N-6} \\
&&&\\
\hline
&&&\\
k_1 +  k_3 \ , \  - k_2 \pm k_i & i=4,... , N-3 &  & 2N-11\\
&&(1,-1)&\\
\frac{1}{2} (k_1 - k_2 + k_3 \pm k_4 \pm ... \pm {\bf u}) & \text{even \# of +} & & 2^{N-6} \\
&&&\\
\hline
\end{array}
\end{equation*}
\vspace{-16pt}
\caption{The Magic Star for $\esn$; ${\bf u}:=k_{N-2}+k_{N-1}+k_N$}\label{t:magices}
\end{table}

\begin{table}[h]
\begin{equation*}
\begin{array}{| c l c | c | c |}
\hline
\hspace{2.9cm} generalized\ roots & & (r,s) & \#\ of\ roots\\
\hline
\pm (k_1 -  k_2) &  & \pm (2,0) & 2\\
\pm (k_2 -  k_3) &  & \pm (-1,3) & 2\\
\pm (k_3 -  k_1) & & \pm (-1,-3) & 2\\
\hline
&&&\\
\pm {\bf v}\ ({\bf v}:=k_{N-1}+k_N)\ , \ \pm k_i \pm  k_j & 4 \le i < j \le N-2 &  & 2(N^2-11N+31)\\
&&(0,0)&\\
\frac{1}{2} (\pm (k_1 + k_2 + k_3) \pm k_4 \pm ... \pm {\bf v}) & \text{even \# of +} & & 2^{N-4} \\
&&&\\
\hline
&&&\\
k_1 +  k_2 \ , \  - k_3 \pm k_i & i=4,... , N-2 &  & 2N-9\\
&&(0,2)&\\
\frac{1}{2} (k_1 + k_2 - k_3 \pm k_4 \pm ... \pm {\bf v}) & \text{even \# of +} &  & 2^{N-5} \\
&&&\\
\hline
&&&\\
- k_1 -  k_2 \ , \  k_3 \pm k_i & i=4,... , N-2 &  & 2N-9\\
&&(0,-2)&\\
\frac{1}{2} (- k_1 - k_2 + k_3 \pm k_4 \pm ... \pm {\bf v}) & \text{even \# of +} & & 2^{N-5} \\
&&&\\
\hline
&&&\\
- k_2 -  k_3 \ , \  k_1 \pm k_i & i=4,... , N-2 &  & 2N-9\\
&&(1,1)&\\
\frac{1}{2} (k_1 - k_2 - k_3 \pm k_4 \pm ... \pm {\bf v}) & \text{even \# of +} &  & 2^{N-5} \\
&&&\\
\hline
&&&\\
k_2 +  k_3 \ , \ - k_1 \pm k_i & i=4,... , N-2 &  & 2N-9\\
&&(-1,-1)&\\
\frac{1}{2} (- k_1 + k_2 + k_3 \pm k_4 \pm ... \pm {\bf v}) & \text{even \# of +} & & 2^{N-5} \\
&&&\\
\hline
&&&\\
- k_1 -  k_3 \ , \  k_2 \pm k_i & i=4,... , N-2 &  & 2N-9\\
&&(-1,1)&\\
\frac{1}{2} (- k_1 + k_2 - k_3 \pm k_4 \pm ... \pm {\bf v}) & \text{even \# of +} &  & 2^{N-5} \\
&&&\\
\hline
&&&\\
k_1 +  k_3 \ , \  - k_2 \pm k_i & i=4,... , N-2 &  & 2N-9\\
&&(1,-1)&\\
\frac{1}{2} (k_1 - k_2 + k_3 \pm k_4 \pm ... \pm {\bf v}) & \text{even \# of +} & & 2^{N-5} \\
&&&\\
\hline
\end{array}
\end{equation*}
\vspace{-16pt}
\caption{The Magic Star for $\estn$; ${\bf v}:=k_{N-1}+k_N$}\label{t:magicest}
\end{table}

\begin{table}[h]
\begin{equation*}
\begin{array}{| c l c | c | c |}
\hline
\hspace{2.9cm} generalized\ roots & & (r,s) & \#\ of\ roots\\
\hline
\pm (k_1 -  k_2) &  & \pm (2,0) & 2\\
\pm (k_2 -  k_3) &  & \pm (-1,3) & 2\\
\pm (k_3 -  k_1) & & \pm (-1,-3) & 2\\
\hline
&&&\\
 \pm k_i \pm  k_j & 4 \le i < j \le N &  & 2(N-3)(N-4)\\
&&(0,0)&\\
\frac{1}{2} (\pm (k_1 + k_2 + k_3) \pm k_4 \pm ... \pm k_N) & \text{even \# of +} & & 2^{N-3} \\
&&&\\
\hline
&&&\\
k_1 +  k_2 \ , \  - k_3 \pm k_i & i=4,... , N &  & 2N-5\\
&&(0,2)&\\
\frac{1}{2} (k_1 + k_2 - k_3 \pm k_4 \pm ... \pm k_N) & \text{even \# of +} &  & 2^{N-4} \\
&&&\\
\hline
&&&\\
- k_1 -  k_2 \ , \  k_3 \pm k_i & i=4,... , N &  & 2N-5\\
&&(0,-2)&\\
\frac{1}{2} (- k_1 - k_2 + k_3 \pm k_4 \pm ... \pm k_N) & \text{even \# of +} & & 2^{N-4} \\
&&&\\
\hline
&&&\\
- k_2 -  k_3 \ , \  k_1 \pm k_i & i=4,... , N &  & 2N-5\\
&&(1,1)&\\
\frac{1}{2} (k_1 - k_2 - k_3 \pm k_4 \pm ... \pm k_N) & \text{even \# of +} &  & 2^{N-4} \\
&&&\\
\hline
&&&\\
k_2 +  k_3 \ , \ - k_1 \pm k_i & i=4,... , N &  & 2N-5\\
&&(-1,-1)&\\
\frac{1}{2} (- k_1 + k_2 + k_3 \pm k_4 \pm ... \pm k_N) & \text{even \# of +} & & 2^{N-4} \\
&&&\\
\hline
&&&\\
- k_1 -  k_3 \ , \  k_2 \pm k_i & i=4,... , N &  & 2N-5\\
&&(-1,1)&\\
\frac{1}{2} (- k_1 + k_2 - k_3 \pm k_4 \pm ... \pm k_N) & \text{even \# of +} &  & 2^{N-4} \\
&&&\\
\hline
&&&\\
k_1 +  k_3 \ , \  - k_2 \pm k_i & i=4,... , N &  & 2N-5\\
&&(1,-1)&\\
\frac{1}{2} (k_1 - k_2 + k_3 \pm k_4 \pm ... \pm k_N) & \text{even \# of +} & & 2^{N-4} \\
&&&\\
\hline
\end{array}
\end{equation*}
\vspace{-16pt}
\caption{The Magic Star for $\eon$}\label{t:magiceo}
\end{table}

\begin{table}[h]
\begin{equation*}
\hspace{-40pt}
\begin{array}{| c l c | c | c |}
\hline
\hspace{2.9cm} generalized\ roots & & (r,s) & \#\ of\ roots\\
\hline
\pm (k_4 -  k_5) &  & \pm (2,0) & 2\\
\pm (k_5 -  k_6) &  & \pm (-1,3) & 2\\
\pm (k_6 -  k_4) & & \pm (-1,-3) & 2\\
\hline
&&&\\
 \pm k_i \pm  k_j & 7 \le i < j \le N &  & 2(N-6)(N-7)\\
&&(0,0)&\\
\frac{1}{2} (\pm (k_1 + k_2 + k_3) \pm (k_4 + k_5 + k_6) \pm k_7... \pm k_N) & \text{even \# of +} & & 2^{N-5} \\
&&&\\
\hline
&&&\\
k_4 +  k_5 \ , \  - k_6 \pm k_i & i=7,... , N &  & 2N-11\\
&&(0,2)&\\
\frac{1}{2} (\pm(k_1 + k_2 + k_3) + k_4 + k_5 -k_6 \pm k_7 ... \pm k_N) & \text{even \# of +} &  & 2^{N-6} \\
&&&\\
\hline
&&&\\
-k_4 -  k_5 \ , \  k_6 \pm k_i & i=7,... , N &  & 2N-11\\
&&(0,-2)&\\
\frac{1}{2} (\pm(k_1 + k_2 + k_3) - k_4 - k_5 + k_6 \pm k_7 ... \pm k_N) & \text{even \# of +} &  & 2^{N-6} \\
&&&\\
\hline&&&\\
-k_5 -  k_6 \ , \  k_4 \pm k_i & i=7,... , N &  & 2N-11\\
&&(1,1)&\\
\frac{1}{2} (\pm(k_1 + k_2 + k_3) + k_4 - k_5 -k_6 \pm k_7 ... \pm k_N) & \text{even \# of +} &  & 2^{N-6} \\
&&&\\
\hline
k_5+  k_6 \ , \  -k_4 \pm k_i & i=7,... , N &  & 2N-11\\
&&(-1,-1)&\\
\frac{1}{2} (\pm(k_1 + k_2 + k_3) - k_4 + k_5 + k_6 \pm k_7 ... \pm k_N) & \text{even \# of +} &  & 2^{N-6} \\
&&&\\
\hline
-k_4 -  k_6 \ , \  k_5 \pm k_i & i=7,... , N &  & 2N-11\\
&&(-1,1)&\\
\frac{1}{2} (\pm(k_1 + k_2 + k_3) - k_4 + k_5 -k_6 \pm k_7 ... \pm k_N) & \text{even \# of +} &  & 2^{N-6} \\
&&&\\
\hline
k_4 +  k_6 \ , \  -k_5 \pm k_i & i=7,... , N &  & 2N-11\\
&&(1,-1)&\\
\frac{1}{2} (\pm(k_1 + k_2 + k_3) + k_4 - k_5 + k_6 \pm k_7 ... \pm k_N) & \text{even \# of +} &  & 2^{N-6} \\
&&&\\
\hline
\end{array}
\end{equation*}
\vspace{-16pt}
\caption{The Magic Star of $\esn$ in the center of the Magic Star of $\eon$}\label{t:magiceseo}
\end{table}

All these sets of generalized roots form a \textquotedblleft Magic Star" as reported in figure \ref{fig:ms}, once projected on the plane spanned by $k_1-k_2$ and $k_1+k_2-2k_3$, with subdivision and grouping of the points and center of the star as reported in tables \ref{t:magicfq}, \ref{t:magices}, \ref{t:magicest}, \ref{t:magiceo}, for $\fqn$, $\esn$, $\estn$ and $\eon$ respectively. In the tables, as well as in figure \ref{fig:ms}, $(r,s)$ denote the pair of scalar products of each root with $k_1-k_2$ and $k_1+k_2-2k_3$, respectively; furthermore, in table \ref{t:magices}  ${\bf u}:=k_{N-2}+k_{N-1}+k_N$, while in in table \ref{t:magicest}  ${\bf v}:=k_{N-1}+k_N$. As a consistency check, it should be remarked here that for $n=1$ one
retrieves the \textquotedblleft Magic Star" projection (or $\mathbf{g}_{2}$
decomposition, as Mukai worded it \cite{Mukai}) of finite dimensional
exceptional Lie algebras \cite{Truini} (cfr. Prop. III.2 below).

By looking at these tables one can readily check that, upon a relabelling of the $k$'s, $\esn$ is the center of the Magic Star of $\eon$ (as reported explicitly in table \ref{t:magiceseo}) and that $\estn=\esn\oplus T_{(r,s)}\oplus T_{(-r,-s)}$, for a fixed pair $(r,s)\in \{(1,1),(-1,1),(0,-2)\}$, where $\esn$ is the center of the Magic Star and $T_{(r,s)}$ is the $(r,s)$ set of roots in table \ref{t:magiceo} of $\eon$. The {\it rank} of the \textquotedblleft Magic Star algebras" $\gdn$, $\fqn$, $\esn$, $\estn$, $\eon$ is defined as the dimension of the vector space $V$ spanned by their roots, namely $2$, $N-4$, $N-2$,  $N-1$, $N$ respectively. By abuse of definition we shall often say {\it root}, for short, instead of
generalized root.\\

From now on, we will restrict to $\esn$, $\estn$, $\eon$ ($\fqn$ deserves a separate treatment, see \cite{EP3}) and denote by $\lms$ anyone of them, by $\Phi$ the set of generalized roots of $\lms$ and by $R$ its rank. We recall that $N=4(n+1)$, $n=1,2,...$, hence $R=N-2=4n+2$ for $\esn$, $R=N-1=4n+3$ for $\estn$, $R=N=4n+4$ for $\eon$.

We denote by $\Phi_O$ and $\Phi_S$ the following subsets of $\Phi$:
\bea{l}
\Phi_O = \{ (\pm k_i \pm  k_j) \in \Phi \}
\\ \\
\Phi_S =  \{ \frac{1}{2} (\pm k_1 \pm k_2 \pm ... \pm k_N)  \in \Phi \}
\eea

\begin{remark}\label{r:dd} Notice that $\Phi_O$ is the root system of $\ddnt$ in the case of $\esn$, of $\ddnd\oplus \au$ in the case of $\estn$ and of $\ddn$ in the case of $\eon$. The corresponding vector spaces of representations in terms of the level $n$ are respectively given by:%
\begin{eqnarray}
\mathbf{e}_{\mathbf{6}}^{(\mathbf{n})} &:&=\overline{\mathbf{\psi }}_{%
\mathbf{d}_{\mathbf{4n+1}}}\oplus \left( \mathbf{d}_{\mathbf{4n+1}}\oplus
\mathbb{C}\right) _{0}\oplus \mathbf{\psi }_{\mathbf{d}_{\mathbf{4n+1}}};
\label{e6n} \\
\mathbf{e}_{\mathbf{7}}^{(\mathbf{n})} &:&=\left( \mathbf{d}_{\mathbf{4n+2}%
}\oplus \mathbf{a}_{1}\right) \oplus \left( \mathbf{\psi }_{\mathbf{d}_{%
\mathbf{4n+2}}},\mathbf{2}\right) =\mathbf{1}_{-2}\oplus \left( \mathbf{\psi
}_{\mathbf{d}_{\mathbf{4n+2}}}\right) _{-1}\oplus \left( \mathbf{d}_{\mathbf{%
4n+2}}\oplus \mathbb{C}\right) _{0}\oplus \left( \mathbf{\psi }_{\mathbf{d}_{%
\mathbf{4n+2}}}\right) _{1}\oplus \mathbf{1}_{2};  \label{e7n} \\
\mathbf{e}_{\mathbf{8}}^{(\mathbf{n})} &:&=\mathbf{d}_{\mathbf{4n+4}}\oplus
\mathbf{\psi }_{\mathbf{d}_{\mathbf{4n+4}}}=\left( \mathbf{6+8n}\right)
_{-2}\oplus \left( \overline{\mathbf{\psi }}_{\mathbf{d}_{\mathbf{4n+3}%
}}\right) _{-1}\oplus \left( \mathbf{d}_{\mathbf{4n+3}}\oplus \mathbb{C}%
\right) _{0}\oplus \left( \mathbf{\psi }_{\mathbf{d}_{\mathbf{4n+3}}}\right)
_{1}\oplus \left( \mathbf{6+8n}\right) _{2},  \nonumber \\
&&  \label{e8n}
\end{eqnarray}%
where $\mathbf{\psi }_{\mathbf{d}_{\mathbf{4n+1}}}\equiv \boldsymbol{2}^{%
\mathbf{4n}}$, $\mathbf{\psi }_{\mathbf{d}_{\mathbf{4n+2}}}\equiv
\boldsymbol{2}^{\mathbf{4n+1}}$, $\mathbf{\psi }_{\mathbf{d}_{\mathbf{4n+3}%
}}\equiv \boldsymbol{2}^{\mathbf{4n+2}}$ and $\mathbf{\psi }_{\mathbf{d}_{%
\mathbf{4n+4}}}\equiv \boldsymbol{2}^{\mathbf{4n+3}}$ respectively denote
the Weyl semispinors of $\mathbf{d}_{\mathbf{4n+1}}$, $\mathbf{d}_{\mathbf{%
4n+2}}$, $\mathbf{d}_{\mathbf{4n+3}}$ and $\mathbf{d}_{\mathbf{4n+4}}$. In
the suitable real cases, (\ref{e6n}), (\ref{e7n}) and (\ref{e8n})
respectively determine the 3-grading of $\mathbf{e}_{\mathbf{6}}^{(\mathbf{n}%
)}$, the 5-grading (of contact type) of $\mathbf{e}_{\mathbf{7}}^{(\mathbf{n}%
)}$, and the 5-grading (of extended Poincar\'{e} type) of $\mathbf{e}_{%
\mathbf{8}}^{(\mathbf{n})}$; in turn, for $n=1$ these reproduce the graded
structure of $\mathbf{e}_{\mathbf{6}}$, $\mathbf{e}_{\mathbf{7}}$, and $%
\mathbf{e}_{\mathbf{8}}$, respectively (cfr. \cite{Group32-1} for further
discussion and details).
\end{remark}

\begin{proposition} For all $\rho\in \Phi_O$ and $x\in \Phi$:
 $2\dfrac{(x, \rho)}{(\rho,\rho)}\in \zz$ and $w_\rho(x) = x - 2\dfrac{(x, \rho)}{(\rho,\rho)}\rho\in \Phi$ (the set of generalized roots is closed under the Weyl reflections by all $\rho\in \Phi_O$).  The set of generalized roots is closed under the Weyl reflections by all $\rho\in \Phi$ if and only if $n=1$.
\end{proposition}
\nin {\bf Proof.} If $\rho\in \Phi_O$ then $(\rho,\rho)=2$ and $(x,\rho)\in \{0,\pm1,\pm 2\}$, hence $2\dfrac{(x, \rho)}{(\rho,\rho)}\in \zz$ and $w_\rho(x) = x - (x, \rho) \rho$. If $(x,\rho)=0$ then $w_\rho(x) = x\in \Phi$. If $(x,\rho)=\pm 1$ then $w_\rho(x) = x\mp \rho \in \Phi$ as we shall prove in proposition \ref{sproots}. If $(x,\rho)=\pm 2$ then $\rho=\pm x$ and $w_\rho(x) = - x\in \Phi$. Suppose now that both $x,\rho\in \Phi_S$ and write $x= \frac12 \sum \lambda_i k_i$, $\rho= \frac12 \sum \mu_i k_i$ where $\lambda_i,\mu_i\in \{-1,1\}$. We can certainly pick an $x$ such that $(x,\rho)=-n$ which occurs whenever for two  indices $j,\ell$ $\lambda_j=\mu_j$ and $\lambda_\ell=\mu_\ell$ while $\lambda_i = - \mu_i$ for $i\ne j,\ell$. Hence $w_\rho(x) = x + 2\dfrac{n}{n+1}\rho = \frac12\sum \nu_i k_i$. We have that $|\nu_j|= \left| \lambda_j + 2\dfrac{n}{n+1} \mu_j \right|= \dfrac{3n+1}{n+1}\ge 2$ and $|\nu_j| = 2$ if and only if $n=1$, in which case $\Phi$ is the root system of a simple Lie algebra. For $n>1$ there is no root with such a coefficient $\nu_j$.\hfill $\square$\\

It should be pointed out that the level $n$ parametrizes the mod.8 Bott periodicity of the Clifford structures (\textit{cfr. e.g.} \cite{Spinor} and Refs. therein) corresponding to the generalized roots $\Phi_S$, which sit into a (semi)spinor representation of the orthogonal Lie algebra whose $\Phi_O$ is the root lattice. Such a mod. 8 periodicity in the framework of the generalization of exceptional Lie algebras provided by persistence of the Magic Star projection, justifies the name \textquotedblleft Exceptional Periodicity" which we adopted since \cite{Mile-High-EP} to describe this mathematical framework.\\

We introduce the basis $\Delta = \{\alpha_1 , ... ,\alpha_R\}$ of $\Phi$, with $\alpha_i = k_i-k_{i+1}\, , \ 1\le i \le R-2$, $\alpha_{R-1} = k_{R-2}+k_{R-1}$ and $\alpha_R = -\frac12(k_1+k_2+...+k_N)$; we order $\Delta$ by setting $\alpha_i < \alpha_{i+1}$:
\be\label{sroots}
\Delta = \{k_1-k_2 < k_2-k_3< ...< k_{R-2} - k_{R-1}< k_{R-2} + k_{R-1}< -\um(k_1+k_2+...+k_N)\}
\ee

\begin{proposition}\label{approp} The set $\Delta$ in \eqref{sroots} is a set of {\bf simple generalized roots}, by which we mean:

\bit
\item[i)] $\Delta$ is a basis of the Euclidean space $V$ of finite dimension $R$;
\item[ii)] every root $\beta$ can be written as a linear combination of roots of $\Delta$ with all positive or all negative integer coefficients: $\beta = \sum \ell_i \alpha_i$ with $\ell_i \ge 0$ or $\ell_i\le 0$ for all $i$.
\eit
\end{proposition}

\nin {\bf Proof.} The set $\Delta = \{\alpha_1 , ... ,\alpha_R\}$ is obviously a basis in $V$.
Let ${\bf u} = k_{N-2}+k_{N-1}+k_N, k_{N-1}+k_N, k_N$ for $\esn$, $\estn$, $\eon$ respectively. We have:
\bea{ccl} \label{kal}
k_{R-1} &= &\frac12 (\alpha_{R-1} - \alpha_{R-2})\\ \\
k_i &= &\alpha_i + k_{i+1} = \sum_{\ell=i}^{R-2}\alpha_\ell +\frac12 (\alpha_{R-1}-\alpha_{R-2}) \ , \ 1\le i\le R-2\\ \\
{\bf u} &= &- 2 \alpha_R - \sum_{\ell=1}^{R-2} {\ell \alpha_\ell - \frac{R-1}{2}(\alpha_{R-1}}-\alpha_{R-2})
\eea
from which we obtain, for $1\le i<j\le R-1$ and forcing $\sum_{\ell=r}^s{\alpha_\ell}=0$ if $r>s$:

\be\label{posroots1}
\left.
\begin{array}{rcl}
k_i + k_j &= & \sum_{\ell=i}^{R-3}\alpha_\ell + \sum_{\ell=j}^{R-2}\alpha_\ell + \alpha_{R-1}\\ \\
k_i - k_j &= & \sum_{\ell=i}^{j-1}\alpha_\ell
\end{array}
\right\}
1\le i<j\le R-1
\ee
for $\estn$, namely for $R=N-1=4n+3$ and ${\bf u}=k_{N-1}+k_N$:
\bea{rcl}\label{posroots2}
-{\bf u} &= & 2 \alpha_R + \sum_{\ell=1}^{R-3} {\ell \alpha_\ell} + 2n\, \alpha_{R-2} + (2n+1)\alpha_{R-1} \ , \text{ for } \estn
\eea
for $\eon$, namely for $R=N$:
\bea{lcl}
\pm k_i - k_N &= & 2 \alpha_N +\sum_{\ell=1}^{i-1} {\ell \alpha_\ell}+ \sum_{\ell=i}^{N-3} {(\ell\pm1)\alpha_\ell} + (2n+ \frac{1 \pm 1}{2})\alpha_{N-2}\\ \\
&&+ (2n+1 + \frac{1 \pm 1}{2})\alpha_{N-1} \ , \quad i\le N-2\\ \\
\pm k_{N-1} - k_N &= & 2 \alpha_N + \sum_{\ell=1}^{N-3} {\ell \alpha_\ell} + (2n+\frac{1 \mp 1}{2})\alpha_{N-2}+ (2n+1 + \frac{1 \pm 1}{2})\alpha_{N-1}
\label{posroots3}
\eea

We see that all the roots in \eqref{posroots1},\eqref{posroots2},\eqref{posroots3} are the sum of simple roots with all positive integer coefficients.
These are half of the roots in $\Phi_O$ and they are all positive roots. \\
The rest of the roots in $\Phi_O$ are negative and they are obviously the sum of simple roots with integer coefficients that are all negative.

Finally all the roots in $\Phi_S$ that contain $-\frac12 {\bf u}$ can be obtained from $\alpha_R$ by flipping an even number of signs and this is done by adding to $\alpha_R$ a certain number of terms of the type $k_i+k_j$, $1\le i<j\le R-1$. These are all positive roots, a half of the roots in $\Phi_S$ and are linear combination of simple roots with integer coefficients that are all positive.\\
The negative roots are similarly obtained by adding to $-\alpha_N$ a certain number of terms of the type $-(k_i+k_j),\, 1\le i<j\le R-1$, and are linear combination of simple roots with integer coefficients that are all negative. \hfill $\square$\\

\begin{remark}\label{mR}
A consequence of the proof of Proposition \ref{approp} is that for all roots $\beta = \sum_{i = 1}^R {m_i\alpha_i}$ the coefficient $m_R$ is such that:
\bea{ll}
m_R \in \{0,\pm 2\} & \text{if } \beta \in \Phi_O\\
m_R =\pm 1 & \text{if } \beta \in \Phi_S\\ \label{rema}
\eea
\end{remark}

\begin{proposition} \label{sproots}
For each $\alpha \in \Phi_O , \beta \in \Phi$ the scalar product $(\alpha,\beta) \in \{\pm 2, \pm 1, 0\}$; $\alpha + \beta$ ( respectively $\alpha - \beta$) is a root if and only if $(\alpha , \beta) = -1$ (respectively $+1$); if both $\alpha+\beta$ and  $\alpha-\beta$ are not in $\Phi\cup\{0\}$ then  $(\alpha,\beta)=0$.\\
 For each $\alpha , \beta \in \Phi_S$ the scalar product $(\alpha,\beta) \in \{\pm (n+1), \pm n , \pm (n-1), ..., 0 \}$; $\alpha + \beta$ ( respectively $\alpha - \beta$) is a root if and only if $(\alpha , \beta) = -n$ (respectively $+n$).\\
For $\alpha , \beta \in \Phi$ if $\alpha+\beta$ is a root then $\alpha - \beta$ is not a root.
\end{proposition}

\nin {\bf Proof.} If both $\alpha , \beta \in \Phi_O$ the proof of the whole Proposition follows from the fact that $\Phi_O \subseteq \Phi_\ddn$, where $\Phi_\ddn$ is the root system of $\ddn$, which is a simply laced Lie algebra. If $\alpha \in \Phi_O , \beta \in \Phi_S$ then $(\alpha,\beta) \in \{\pm 1, 0\}$, as a trivial computation explicitly shows. Moreover, let us write $\alpha = \sigma_i k_i +\sigma_j k_j \, , \ i < j$, $\sigma_{i,j}\in \{-1,1\}$. Then $(\alpha , \beta) = -1$ if and only if $\beta = \frac12 (\pm k_1 \pm ... -\sigma_i k_i \pm ... -\sigma_j k_j \pm ... \pm k_N)$, which is true if and only if $\alpha + \beta \in \Phi$ (in particular $\alpha + \beta \in \Phi_S$). Similarly $(\alpha , \beta) = 1$ if and only if $\alpha - \beta \in \Phi$. As a consequence, if both $\alpha\pm\beta \notin \Phi\cup\{0\}$ then  $(\alpha,\beta)\ne \pm 1$ and also $(\alpha,\beta)\ne \pm 2$ because $(\alpha,\beta)= \pm 2$ if and only if $\alpha = \pm \beta$; therefore $(\alpha,\beta)=0$.\\
If both $\alpha , \beta \in \Phi_S$, then all their signs but an even number $2m$ must be equal, $m=0,...,N/2=2(n+1)$ and we get $(\alpha,\beta) = \frac14 (N-2m - 2m) = n+1-m = n+1, n, n-1, ... ,-(n+1)$ for  $m=0,...,2(n+1)$. Moreover, since $\pm k_i \pm k_i  \in \{0,\pm 2 k_i\} \, , \ i = 1, ... N$ then $\alpha + \beta \in \Phi$ if and only if all signs are opposite but 2 (in which case $\alpha + \beta$ is actually in $\Phi_O$) and this is true if and only if $(\alpha , \beta) = -\frac14(N-4)=-n$. Similarly $(\alpha , \beta) = n$ if and only if $\alpha - \beta \in \Phi$. The last statement of the Proposition follows trivially. \hfill $\square$\\


\section{The Magic Star algebra $\lms$}\label{s:lms}

We define the Magic Star algebra $\lms$ (as before $\lms$ is either $\esn$ or $\estn$ or $\eon$) by extending the construction of a Lie algebra from a root system, \cite{carter} \cite{hum} \cite{graaf}. In particular, we generalize the algorithm in \cite{graaf} for simply laced Lie algebras, since also in our set of generalized roots the $\beta$ chain through $\alpha$, namely the set of roots $\alpha+c\beta$, $c\in \zz$, has length one.\\

We give $\lms$ an algebra structure of rank $R$ over a field extension $\fff$ of the rational integers $\zz$ in the following way\footnote{Specifically, we will take $\fff$ to be the complex field $\cc$.}:
\bit
\item[a)]  we select the set of simple generalized roots $\Delta = \{\alpha_1 , ... ,\alpha_R\}$ of $\Phi$
\item[b)] we select a basis $\{ h_1 ,...,h_R\}$ of the $R$-dimensional vector space $H$ over $\fff$ and set $h_\alpha = \sum_{i=1}^R c_i h_i$ for each $\alpha  \in \Phi$ such that $\alpha = \sum_{i=1}^R c_i \alpha_i$
\item[c)] we associate to each $\alpha  \in \Phi$ a one-dimensional vector space $L_\alpha$ over $\fff$ spanned by $x_\alpha$
\item[d)] we define $\lms = H \bigoplus_{\alpha \in \Phi} {L_\alpha}$ as a vector space over $\fff$
\item[e)] we give $\lms$ an algebraic structure by defining the following multiplication on the basis
$\blms = \{ h_1 ,...,h_R\} \cup \{x_\alpha \ | \ \alpha \in \Phi\}$, extended by linearity to a bilinear multiplication $\lms\times \lms\to \lms$:
	\be\begin{array}{ll}
	&[h_i,h_j] = 0 \ , \ 1\le i, j  \le R \\
	&[h_i , x_\alpha] = - [x_\alpha , h_i] = (\alpha, \alpha_i )\, x_\alpha \ , \ 1\le i \le R \ , \ \alpha \in \Phi \\
	&[x_\alpha, x_{-\alpha} ] = - h_\alpha\\
	&[x_\alpha,x_\beta] = 0 \ \text{for } \alpha, \beta \in \Phi \ \text{such that } \alpha + \beta \notin			 	\Phi \ \text{and } \alpha \ne - \beta\\
	&[x_\alpha,x_\beta] = \varepsilon (\alpha , \beta)\,  x_{\alpha+\beta}\ \text{for } \alpha , \beta \in \Phi \ \text{such that }  \alpha+ \beta \in \Phi\\
	\end{array} \label{comrel}\ee
\eit

where $\varepsilon (\alpha , \beta)$ is the {\it asymmetry function}, introduced in \cite{kac}, see also \cite{graaf}, defined as follows:\\
\begin{definition} Let $\Ll$ denote the lattice of all linear combinations of the simple generalized roots with integer coefficients
\be
\Ll = \left\{ \sum_{i=1}^R c_i \alpha_i \ |\ c_i \in \zz \ , \ \alpha_i \in \Delta \right\}
\label{lattice}
\ee
the asymmetry function $\varepsilon (\alpha , \beta): \ \Ll \times \Ll \to \{-1,1\}$ is defined by:
\be\label{epsdef}
\varepsilon (\alpha , \beta) = \prod_{i,j=1}^R \varepsilon (\alpha_i , \alpha_j)^{\ell_i m_j} \quad \text{for } \alpha = \sum_{i=1}^R \ell_i\alpha_i \ ,\ \beta = \sum_{j=1}^R m_j \alpha_j
\ee
where $\alpha_i , \alpha_j \in \Delta$ and
\be\label{epsdef1}
\varepsilon (\alpha_i , \alpha_j) = \left\{
\begin{array}{ll}
-1 & \text{if } i=j\\ \\
-1 & \text{if } \alpha_i + \alpha_j  \text{ is a root and } \alpha_i < \alpha_j\\ \\
+ 1 & \text{otherwise}
\end{array}
\right.
\ee
\end{definition}

Note that (\ref{e6n})-(\ref{e8n}) exhibit Bott periodicity (due to the
increasing mod. $4$ of the rank of the corresponding lattice, or
equivalently to the increasing mod.$8$ of the argument of the corresponding $%
\mathbf{d}$-type Lie algebra). The commutation relations of the
corresponding generators are given in terms of the asymmetry function
defined in Definition 4.1. It is here worth anticipating that Magic Star
algebras $\mathfrak{L}_{MS}$ are not simply non-reductive, spinorial
extensions of Lie algebras, but rather they are characterized by a
non-translational (i.e., non-Abelian) nature of their spinorial sector; this
implies that they are Lie algebras only for $n=1$, i.e. at the trivial level
of Exceptional Periodicity, whereas for $n\geqslant 2$ they are not Lie
algebras, because the Jacobi identity is violated in the spinorial sector
itself (for this, we address the reader to the discussion in \cite{EP2}).

\section{Properties of the asymmetry function}\label{s:epsilon}

We now show some properties of the asymmetry function $\varepsilon (\alpha , \beta): \ \Ll \times \Ll \to \{-1,1\}$, that is crucial in the definition of the algebra $\lms$. In particular we show that, for $\as,\bb,\as+\bb\in\Phi$, $\varepsilon (\alpha , \beta) = - \varepsilon (\beta , \alpha)$ which implies that the bilinear product \eqref{comrel} is antisymmetric.

\begin{proposition} \label{epsprop}The asymmetry function $\varepsilon$ satisfies, for $\alpha , \beta, \gamma , \delta \in \Ll$, $\alpha = \sum{m_i \alpha_i}$ and $\beta = \sum{n_i \alpha_i}$:
\bes\begin{array}{rrcl}
i) & \varepsilon (\alpha + \beta, \gamma) & =& \varepsilon (\alpha , \gamma)\varepsilon (\beta , \gamma) \\
ii) &\varepsilon (\alpha ,\gamma + \delta) & = &\varepsilon (\alpha , \gamma)\varepsilon (\alpha , \delta) \\
iii) &\varepsilon (\alpha , \alpha) & = &(-1)^{\frac12 (\alpha ,\alpha)-m_R^2 \frac{n-1}2}\\
iv) &\varepsilon (\alpha , \beta) \varepsilon (\beta , \alpha) &=& (-1)^{(\alpha , \beta) - m_R n_R (n-1)}\\
v) &\varepsilon (0 , \beta) &=& \varepsilon (\alpha , 0) = 1 \\
vi) &\varepsilon (-\alpha , \beta) &=& \varepsilon (\alpha , \beta)^{-1} = \varepsilon (\alpha , \beta) \\
vii) &\varepsilon (\alpha , -\beta) &=& \varepsilon (\alpha , \beta)^{-1} = \varepsilon (\alpha , \beta) \\
\end{array}\ees
\end{proposition}

\nin {\bf Proof.} The first two properties follow directly from the definition. In order to prove $iii)$ we first notice that for $\alpha_i,\alpha_j\in \Delta$ , $i\ne j$, $(\alpha_i, \alpha_j) \in \{0,-1\}$. Therefore:
\be\begin{array}{rcl}
\varepsilon (\alpha , \alpha) &=& \prod_{1\le i,j\le R} \varepsilon (\alpha_i , \alpha_j)^{m_im_j}  =  \prod_{1\le i<j\le R} (-1)^{m_im_j (\alpha_i , \alpha_j)}  \prod_{1\le i \le R} (-1)^{m_i^2}\\ \\
&=& (-1)^{\sum_{1\le i<j\le R}{m_im_j (\alpha_i , \alpha_j)} + \frac12 \sum_{1\le i\le R}{m_i^2(\alpha_i , \alpha_i)} - \frac12  \sum_{1\le i\le R}{m_i^2 ((\alpha_i , \alpha_i)} - 2)} \\ \\
&=& (-1)^{\frac12 (\alpha , \alpha) - \frac12 m_R^2 ((\alpha_R , \alpha_R) - 2)} =  (-1)^{\frac12 (\alpha , \alpha) - m_R^2 \frac{n-1}2}
\end{array}\ee

Property $iv)$ follows from $iii)$ by replacing $\alpha$ with $\alpha+\beta$ and using the first two. If $\alpha = \sum m_i  \alpha_i$, $\beta = \sum n_i \alpha_i$ and $\alpha +\beta = \sum \ell_i \alpha_i = \sum (m_i+n_i) \alpha_i$ we get:

\bea{rcl}
\varepsilon (\alpha + \beta , \alpha +\beta) &=& (-1)^{\frac12 (\alpha +\beta, \alpha+\beta) - \ell_R^2 \frac{n-1}2} =
(-1)^{\frac12 (\alpha, \alpha)  +  \frac12 (\beta, \beta) +  (\alpha, \beta)- \ell_R^2 \frac{n-1}2} \\ \\
&=& (-1)^{\frac12 (\alpha , \alpha) - m_R^2 \frac{n-1}2} (-1)^{\frac12 (\beta, \beta) - n_R^2 \frac{n-1}2}
\varepsilon (\alpha, \beta)\varepsilon (\beta , \alpha )
\eea
from which the property follows. Property $v)$ is a trivial consequence of the definition, whereas properties $vi)$ and $vii)$ follow from property $v)$ together with $i)$ and $ii)$. \hfill$\square$\\

\begin{proposition} If $\alpha, \beta ,\alpha+\beta \in \Phi$ then:
\beas{rll}
i) &\varepsilon (\alpha , \alpha) = -1 & \alpha \in \Phi\\
ii) &\varepsilon (\alpha , \beta) = - \varepsilon (\beta , \alpha) & \alpha,\beta,(\alpha+\beta)\in \Phi\qquad \text{antisymmetry}\\
iii) &\ep(\as ,\bb) = \ep(\bb, \as + \bb) & \text{if } \as, \as+\bb \in \Phi\, , \ \bb\in \Ll \\
iv) &\ep(\as ,\bb) = \ep(\bb, \as - \bb) & \text{if } \as, \as-\bb \in \Phi\, , \ \bb\in \Ll \\
\eeas
\label{remas}
\end{proposition}

\nin {\bf Proof.}
By Remark \ref{mR} if $\alpha\in \Phi_O$ then $(\alpha,\alpha) = 2$ and $m_R^2/2$ is even, hence $\varepsilon (\alpha,\alpha) = -1$. If $\alpha\in \Phi_S$ then $(\alpha,\alpha) = n+1$ and $m_R^2=1$. Therefore if $(-1)^{\frac12 (\alpha ,\alpha)-m_R^2 \frac{n-1}2}= (-1)^{\frac12 (n+1-n+1)} = -1$. As a consequence, if $\alpha, \beta ,\alpha+\beta \in \Phi$, then $-1=\varepsilon (\alpha+\beta ,\alpha+ \beta)=\varepsilon (\alpha , \beta)\varepsilon (\beta,\alpha)$, hence $\varepsilon (\alpha , \beta) = - \varepsilon (\beta , \alpha)$ \\
Finally we prove $iii)$ - a similar proof holds for $iv)$ -:
\bes
\ep(\as ,\bb) = \ep(\as ,\as-\as+\bb) =\ep(\as ,\as)\ep(\as ,\as+\bb) =  - \ep(\as+\bb -\bb ,\as+\bb) = \ep(\bb, \as + \bb)
\ees
\hfill $\square$


\section{Derivations and automorphisms of the Magic Star algebra $\lms$}\label{s:da}

We will now begin to study the properties of the novel mathematical entity given by the Magic Star algebra $\lms$ introduced in previous Sections. To this end, in the present Section we prove that the inner derivations
(or, at group level, the inner
 automorphisms) of the Magic Star algebra $\lms$ are, in the simply laced case under consideration, given by its orthogonal Lie subalgebra.\\

Let us denote by $\dd$ the $\lms$ Lie subalgebra $\ddnt$ if $\lms=\esn$, $\ddnd\oplus\au$ if $\lms=\estn$, $\ddn$ if $\lms=\eon$, spanned by $\{x_\alpha \, , \ \alpha \in \Phi_O\}$; see Remark \sref{r:dd}.

We know that for $n=1$, namely $\lms=\es,\est,\eo$, the adjoint action $\adj_x:y \to [x,y]$ is a derivation of the algebra $\lms$ (and hence $\exp(\zeta \adj_x)$ is an automorphism of $\lms$).  This is due to the Jacobi identity. We have the following result in the case $n>1$:

\begin{proposition} For $n>1$ the adjoint action $\adj_x:y \to [x,y]$ is a derivation of the algebra $\lms$ (and hence $\exp(\zeta \adj_x)$ is an automorphism of $\lms$) if and only if $x\in \dd$.
\label{der}
\end{proposition}
\nin {\bf Proof.} By the linearity of the adjoint action it is sufficient to prove the proposition for basis elements $X_0$ of the algebra $\dd$ and basis elements $X_1, X_2$ of $\lms$, and it amounts to the Jacobi identity:
\be\begin{array}{rl}
J_0 + J_1 + J_2 = 0 &, \ J_\ell := [[X_\ell, X_{\ell +1}], X_{\ell +2}] \  \ \text{(indices mod } 3) \\
 \label{jacobi} \\
X_\ell \in \blms &,\ \ell = 0,1,2 \ , \ X_0\in \dd
\end{array}
\ee

 The identity is trivial if all the $X_\ell$'s are in $H$. If 2 of them are in $H$ and 1 is not it amounts to
 $[[h_i , x_\alpha] , h_j] = [[h_j , x_\alpha] , h_i]$, which holds since both members are equal to $-(\alpha , \alpha_i)(\alpha , \alpha_j)$. If only one of the $X_\ell$'s is in $H$ and 2 are not then:
 \bit
 \item[a1)] if $\alpha+\beta\in \Phi$ the Jacobi identity is equivalent to the identity:\\
  $\varepsilon (\alpha , \beta) (\alpha , \alpha_i) - (\alpha+\beta, \alpha_i) \varepsilon (\alpha , \beta) + (\beta,\alpha_i) \varepsilon (\alpha , \beta) = 0$;
 \item[a2)] if $\alpha = - \beta$ the Jacobi identity is equivalent to $(\alpha , \alpha_i) - (\alpha , \alpha_i) = 0$;
 \item[a3)] if $\alpha+\beta\notin \Phi\cup\{ 0\}$ then the Jacobi identity is trivially satisfied.
 \eit

 \nin From now on none of the $X_\ell$'s is in $H$. Set $X_0 = x_\alpha$, $\alpha\in\Phi_O$, $X_1 = x_\beta$, $X_2 = x_\gamma$ so that $J_0=[[x_\as,x_\bb],x_\gh]$, $J_1=[[x_\bb,x_\gh],x_\as]$, $J_2=
 [[x_\gh,x_\as],x_\bb]$.\\

If none of the sums $\alpha+\beta$, $\alpha+\gamma$, $\beta+\gamma$ is in $\Phi$ nor is 0, then $J_0 +J_1 + J_2 = 0$. Also if any two roots are equal then $J_0 +J_1 + J_2 = 0$. In fact $[[x_\alpha,x_\beta],x_\alpha]+[[x_\beta,x_\alpha],x_\alpha]+
 [[x_\alpha,x_\alpha],x_\beta] = [[x_\alpha,x_\beta],x_\alpha]- [[x_\alpha,x_\beta],x_\alpha]=0$\\

From now on at least one of $\alpha+\beta$, $\alpha+\gamma$, $\beta+\gamma$ is in $\Phi\cup \{0\}$.
Suppose first that $\alpha+\beta=0$. Then $J_0 = - (\gamma, \alpha)x_\gamma$ and we have the following possibilities:
 \bit
 \item[b1)] if $\alpha + \gamma = 0$ or $\alpha - \gamma = 0$ then either $\beta = \gamma$ or $\alpha = \gamma$ in which case the Jacobi identity becomes trivial;
 \item[b2)] if $\alpha + \gamma \in \Phi$ then $\gamma-\beta=\gamma+\alpha\in \Phi$, therefore $\gamma+\beta \notin\Phi \cup \{ 0\}$ and $J_1=0$. We get $J_2 = \varepsilon(\gamma,\alpha)\varepsilon(\alpha+\gamma,-\alpha) x_\gamma = \varepsilon(\alpha,-\alpha) x_\gamma=-x_\gamma$ (because of Propositions \sref{epsprop} and \sref{remas}) and $J_0 = -(\gamma, \alpha)x_\gamma = x_\gamma$ hence Jacobi is verified;
 \item[b3)] if $\alpha - \gamma \in \Phi$ the proof is similar to b2);
 \item[b4)] if both $\alpha \pm \gamma \notin \Phi\cup \{ 0\}$ then $(\alpha,\gamma)=0$, by Proposition \sref{sproots} being $\alpha\in\Phi_O$, hence $J_0 = J_1 = J_2 = 0$.
 \eit
Similarly if we suppose $\alpha+\gamma=0$.\\
Suppose now $\beta+\gamma=0$. Then $J_1 = -(\alpha,\beta)x_\alpha$ and we have the following possibilities:
 \bit
 \item[c1)] if $\alpha + \beta = 0$ or $\alpha - \beta = 0$ then either $\alpha=\beta$ or $\alpha = \gamma$ in which case the Jacobi identity becomes trivial;
 \item[c2)] if $\alpha + \beta \in \Phi$ then $\alpha-\beta=\alpha+\gamma \notin\Phi \cup \{ 0\}$ and $J_2=0$. We get $J_0 = \varepsilon(\alpha,\beta)\varepsilon(\alpha+\beta,\gh) x_\alpha = \varepsilon(\beta,-\beta) x_\alpha=-x_\alpha$ (because of Propositions \sref{epsprop} and \sref{remas}) and $J_1 = -(\alpha,\beta)x_\alpha = x_\alpha$, since $\as\in\Phi_O$, and Jacobi is verified;
 \item[c3)] if $\alpha - \beta \in \Phi$ the proof is similar to c2);
 \item[c4)] if both $\alpha \pm \beta \notin \Phi\cup \{ 0\}$ then $(\alpha,\beta)=0$, by Proposition \sref{sroots} being $\alpha\in\Phi_O$, hence $J_0 = J_1 = J_2 = 0$.
 \eit

 From now on $\alpha,\beta,\gamma\in \Phi$ and $\alpha\ne\pm\beta\ne\pm\gamma\ne\pm\alpha$.\\

If $\alpha+\beta+\gamma\notin \Phi \cup\{ 0 \}$ then $J_0 = J_1 = J_2 = 0$. Consider $[[x_\alpha,x_\beta],x_\gamma]$. If $[x_\alpha,x_\beta]=0$ the statement is proven, otherwise $[x_\alpha,x_\beta] = \varepsilon(\alpha,\beta)x_{\alpha+\beta}$.  Now $[x_{\alpha+\beta},x_\gamma]=0$ because $(\alpha+\beta)+\gamma\notin\Phi$ and the statement is proven.\\

From now on $\alpha,\beta,\gamma\in \Phi$, $\alpha\ne\pm\beta\ne\pm\gamma\ne\pm\alpha$ and $\alpha+\beta+\gamma\in \Phi \cup\{ 0 \}$.

\bit
\item[d1)] Take first $\alpha+\beta+\gamma =  0$. Since $-\alpha\in \Phi$, $-\alpha=\beta+\gamma \in \Phi$ and similarly for $-\beta$, $-\gamma$. Therefore $(\alpha+\beta) , (\beta+\gamma), (\alpha+\gamma)\in \Phi$ and $J_0 = - \eab (h_\alpha +h_\beta)$, $J_1 = - \ebg (h_\beta+h_\gamma) = \ebg h_\alpha$,
 $J_2 = - \ega (h_\alpha+h_\gamma) = \ega h_\beta$. By Proposition \sref{remas} and $\alpha+\beta+\gamma=0$ we get $\eab= \varepsilon(-\beta-\gamma,\beta) = -\varepsilon(\gamma,\beta) = \ebg$ and also $\eab= \varepsilon(\alpha,-\alpha-\gamma) =\ega$ and $J_0 + J_1 + J_2 = 0$.\\
\eit

From now on
\bea{c}
\alpha,\beta,\gamma,\alpha+\beta+\gamma\in \Phi\ ; \ \alpha\ne\pm\beta\ne\pm\gamma\ne\pm\alpha\\
\text{and at least one of } \as+\bb,\bb+\gh, \as+\gh \text{ is in } \Phi
\label{p5}
\eea

 \bit
\item[d2)] By Proposition \sref{only1} the conditions \eqref{p5} imply that exactly two of $\as+\bb,\bb+\gh, \as+\gh$ are in $\Phi$. Suppose $\as+\bb,\bb+\gh \in \Phi$ then\\
\bea{l}
[[x_\alpha,x_\beta],x_\gamma]+[[x_\beta,x_\gamma],x_\alpha] =\\ \eab\eag\ebg+\ebg\eba\ega=\\ \eab\ebg(\eag-\ega)\label{fjac}
\eea
By $iv)$ of Proposition \sref{epsprop} plus Proposition \sref{sproots} and Remark \sref{mR} we have $\eag \ega = (-1)^{(\as,\gh)}$ if $\as\in\Phi_O$. We now prove that $(\as,\gh)=0$ hence $\eag=\ega$. Since $\as+\gh\notin\Phi$ and $\as\in \Phi_O$ we have $(\as,\gh)\in\{0,1\}$. Suppose $(\as,\gh)=1$. Then $(\as+\bb,\gh)\in \{-1,-n\}$ and also $(\as+\bb,\gh)=(\as,\gh)+(\bb,\gh)=1+(\bb,\gh)$. If $\bb$ or $\gh$ are in $\Phi_O$ then
$(\bb,\gh)=-1$ and we get a contradiction. If both $\bb,\gh\in\Phi_S$  then $(\bb,\gh)=-n$ and also $(\as+\bb,\gh)=-n$, again a contradiction. So $(\as,\gh)=0$ and \eqref{fjac} is zero.\\
Similarly if $\as+\gh,\bb+\gh \in \Phi$.\\
Suppose $\as+\bb,\as+\gh \in \Phi$ then\\
\bea{l}
[[x_\as,x_\bb],x_\gh] + [[x_\gh,x_\as],x_\bb]=\\ -\eab\ega\ebg+\ega\egb\eab=\\ \eab\ega(\egb-\ebg)\label{fjacp}
\eea
By $iv)$ of Proposition \sref{epsprop} plus Proposition \sref{sproots} and Remark \sref{mR} we have $\ebg \egb = (-1)^{(\bb,\gh)-m_Nn_N(n-1)}$. If $\bb\in\Phi_O$ or $\gh\in \Phi_O$ then  $\ebg \egb = (-1)^{(\bb,\gh)}$ and $(\bb,\gh)\in\{0,1\}$. Suppose $(\bb,\gh)=1$. Then either $\as+\gh\in\Phi_O$ or $\bb\in\Phi_O$ and $(\as+\gh,\bb)=-1$; but also $(\as+\gh,\bb)=(\as,\bb)+(\gh,\bb)=-1+1=0$, a contradiction. If both $\bb,\gh\in\Phi_S$  then $(\as+\gh,\bb)=-n$ and $(\as,\bb)=-1$. So $(\bb,\gh)=-(n-1)$ and $\ebg \egb = (-1)^{(\bb,\gh)-m_Rn_R(n-1)}= (-1)^{-(n-1)\pm(n-1)}=1$, therefore $\ebg=\egb$ and \eqref{fjacp} is zero.
\eit

We now show that for any $x_\alpha\in \Phi_S$ there exist $x_\bb,x_\gh\in \Phi_S$ such that $[[x_\alpha,x_\beta],x_\gamma]+[[x_\beta,x_\gamma],x_\alpha]+
 [[x_\gamma,x_\alpha],x_\beta]\ne 0$.\\
Let $\as = \um\sum \lambda_i k_i$, $\lambda_i\in\{-1,1\}$, and for a fixed set of different indices $\{j,\ell, m, r,s,t\}$ let $\bb = \lambda_j k_j + \lambda_\ell k_\ell- \as$, $\gh = -\lambda_jk_j - \lambda_\ell k_\ell-\lambda_mk_m - \lambda_rk_r-\lambda_s k_s - \lambda_tk_t+\as$;  then
$\as,\bb,\gh\in\Phi$,  $\as+\bb\in\Phi$,  $\as+\bb+\gh=\lambda_j k_j + \lambda_\ell k_\ell+\gh\in \Phi$ but $\bb+\gh,\gh+\as\notin\Phi$ unless $N=8$, namely $n=1$ (in which case $\gh+\as\in\Phi$) which is excluded in the hypothesis. We thus have $[[x_\alpha,x_\beta],x_\gamma]+[[x_\beta,x_\gamma],x_\alpha]+
 [[x_\gamma,x_\alpha],x_\beta]=[[x_\alpha,x_\beta],x_\gamma]=\eab\varepsilon(\as+\bb,\gh)\ne 0$.\\

Finally, the fact that $\exp(\adj_x)$ is an automorphism if $\adj_x$ is a derivation is a classical result that we recall here.\\
First of all we notice that $\adj_x$ is nilpotent. Let $\der$ be a nilpotent derivation of $\lms$: $\der^r=0$ for some $r$. Then
\be \exp\der = 1+\der+\dfrac{\der^2}{2!}+ ... +\dfrac{\der^{r-1}}{(r-1)!}\ee
 and
\be
\der [x,y]=[\der x,y]+[x,\der y] \, , \quad x,y\in\lms
\ee
imply:
\be
\dfrac1{s!}\der^s [x,y]=
\dfrac1{s!}\sum_{i=0}^s{\binom{s}{i} [\der^i x,\der^{s-i}y]}=\sum_{i=0}^s{\left[\dfrac{\der^i x}{i!},\dfrac{\der^{s-i}y}{(s-i)!}\right]} =\sum_{\substack{i,j\\i+j=s}}{\left[\dfrac{\der^i x}{i!},\dfrac{\der^jy}{j!}\right]}
\ee
hence\\
\bea{rl}
\exp \der [x,y] &=\ \sum_{s\ge0}\sum_{\substack{i,j\\i+j=s}} {\left[\dfrac{\der^i x}{i!},\dfrac{\der^jy}{j!}\right]}
= \sum_{i\ge0}\sum_{j\ge0}{\left[\dfrac{\der^i x}{i!},\dfrac{\der^jy}{j!}\right]}
\\ \\&=\ [\exp \der x,\exp \der y]
\eea
We have used the fact that $\delta^t=0$ for $t\ge r$ implies $\sum_{\substack{i,j}}{\left[\dfrac{\der^i x}{i!},\dfrac{\der^jy}{j!}\right]}=0$ if $i+j\ge r$.\\

This ends the proof of Proposition \sref{der}. \hfill $\square$

\begin{proposition}\label{only1}
Let $\as,\bb,\gh \in \Phi$, $\as\ne\pm \bb\ne\pm \gh\ne\pm\as$, one of which is in $\Phi_O$ and let $\as+\bb\in\Phi$ and $\as+\bb+\gh\in\Phi$. Then one and only one of $(\bb+\gh),(\as+\gh)$ must be in $\Phi$.
\end{proposition}
\proof We will make extensive use of Proposition \sref{sproots}.\\ Suppose both $(\bb+\gh),(\as+\gh)$ be not in $\Phi$.\\
If $\gamma\in \Phi_O$ then $(\bb,\gh),(\as,\gh)\in \{0,1\}$, but $(\as+\bb,\gh)=-1=(\as,\gh)+(\bb,\gh)$, which is impossible.\\
If $\as,\bb\in \Phi_O$ and $\gamma\in \Phi_S$ then, as before, $(\bb,\gh),(\as,\gh)\in \{0,1\}$, but $(\as+\bb,\gh)=-1=(\as,\gh)+(\bb,\gh)$, which is impossible.\\
If $\as\in \Phi_O$ and $\bb,\gamma\in \Phi_S$ then $(\as+\bb)\in \Phi_S$ and $(\as+\bb,\gh)=-n=(\as,\gh)+(\bb,\gh)$, $(\as,\gh)\in\{ 0,1\}$. But $(\as,\gh)=0$ implies $\bb+\gh\in \Phi$ and $(\as,\gh)=1$ implies $(\bb,\gh)=-(n+1)$ hence $\bb=-\gh$, which both contradict the hypothesis.\\
Similarly if $\bb\in \Phi_O$ and $\as,\gh\in \Phi_S$.\\
Finally suppose that both $(\bb+\gh),(\as+\gh)\in \Phi$ and denote by
$s_0, s_1, s_2$ the scalar products $(\alpha+\beta,\gamma)$, $(\alpha,\gamma)$, $(\beta,\gamma)$, respectively.
Since by hypothesis these are all scalar products of roots whose sum is a root, we have by Proposition \sref{sproots}: $s_i\in \{-1,-n\}$, $i=0,1,2$ and $s_0=s_1+s_2\in \{-2,-n-1,-2n\}$.
This is only possible for $n=2$ and $s_0=-2, s_1=s_2=-1$ which implies $\alpha+\beta,\gamma\in \Phi_S$.
But $(\alpha,\gamma)= (\beta,\gamma)=-1$, $(\beta+\gamma), (\alpha+\gamma)\in \Phi$ and $\gamma\in \Phi_S$ imply $\alpha, \beta \in \Phi_O$ hence $(\alpha+\beta)\in \Phi_O$, a contradiction.\\
Viceversa it is always possible to find three roots satisfying the hypothesis of the Proposition such that only two sums of two of them are roots.
\hfill $\square$\\

\section{Further Developments}\label{s:outlook}

The non-Lie, countably infinite chains of finite dimensional generalizations of the exceptional Lie algebras provided by Magic Star algebras pave the way to a number of interesting mathematical developments. Below, we list some of the ones which we plan to report on in the near future.

In the forthcoming papers \cite{EP2}-\nocite{EP3}\cite{EP4}, we will analyze the algebraic structures of the star-shaped projection of Magic Star algebras; remarkably, such structures turn out to be the Hermitian part of the rank-3 matrix algebras introduced by Vinberg in \cite{Vinberg}. Therefore, Exceptional Periodicity not only generalizes exceptional Lie algebras, but also cubic Jordan algebras (and in particular the Albert algebra). Then, we will consider the gradings of Magic Star algebras and the corresponding algebraic structures, which in turn generalize Jordan pairs and Freudenthal triple systems. We will also analyze the non-Lie nature of Magic Star algebras, in particular the subsectors of such algebras which violate the Jacobi identity; we anticipate that such a violation occurs due to the non-trivial (i.e. non-Abelian) nature of the spinorial subsector of the Magic Star algebras.

An interesting line of research stemming from Exceptional Periodicity pertains to study the higher dimensional
weight vectors of algebras akin to lattice vertex algebras (the original
motivation for Borcherds' definition of vertex algebras), that project to a star-shaped, Bott-periodic Magic Star structure. As we have seen, Magic Star algebras are crucially defined by the so-called \textit{asymmetry function}, which acts like the cocycle of a lattice vertex
algebra which gives a twisted group ring $\mathbb{C}_{\epsilon }[\Lambda ]$
over an even lattice $\Lambda $. Correspondingly, the Magic Star algebras span higher-dimensional lattices, beyond that of the self-dual $D=8$ lattice of $\eo$, and thus they potentially allow to probe the symmetries of the heterotic string and moonshine, as well as to formulate a matrix algebra generalization of noncommutative geometry. Remarkably, Exceptional Periodicity provides a novel algebraic method for studying even self-dual
lattices, such as the $\eo\oplus \eo$ and Leech lattices, which already have a well known connection to the Monster vertex algebra and $D=24$
bosonic string compactifications.

Of course, there are also several topics that we are planning to develop in the future, which are strictly related to Quantum Gravity. In particular, a model for interactions based on Exceptional Periodicity which includes gravity and the expansion of space-time. We aim at a new perspective of
elementary particle physics at the early stages of the Universe based on the idea that interactions, defined in a purely algebraic way, are the fundamental objects of the theory, whereas space-time, hence gravity, are derived structures. With an infinite family of new algebras that extend the exceptional Lie algebras, Exceptional Periodicity and Magic Star algebras give a fresh new toolkit for studying emergent spacetime and Quantum Gravity, in dimensions beyond those previously explored, using spectral techniques applied to an infinite class of cubic Hermitian matrix algebras.



\begin{thebibliography}{99}

\bibitem[ABDHN]{ICL-Magic} A. Anastasiou, L. Borsten, M.J. Duff, L.J.
Hughes, S. Nagy, \textit{Super Yang-Mills, division algebras and triality},
JHEP \textbf{1408} (2014) 080, \texttt{arXiv:1309.0546 [hep-th]}. A.
Anastasiou, L. Borsten, M.J. Duff, L.J. Hughes, S. Nagy, \textit{A magic
pyramid of supergravities}, JHEP \textbf{1404} (2014) 178, \texttt{%
arXiv:1312.6523 [hep-th]}.

\bibitem[ABDMN]{ICL-Magic-2} A. Anastasiou, L. Borsten, M.J. Duff, A.
Marrani, S. Nagy, \textit{The Mile High Magic Pyramid}, Contemp. Math.
\textbf{721} (2019) 1-27, \texttt{arXiv:1711.08476 [hep-th]}.

\bibitem[AST]{21} T. Asakawa, S. Sugimoto and S. Terashima, $\mathit{D}$%
\textit{-branes, Matrix Theory and }$\mathit{K}$\textit{-homology}, JHEP
\textbf{0203}, 034 (2002), \texttt{hep-th/0108085}.

\bibitem[AZ]{16} M. Artin and J. J. Zhang, \textit{Noncommutative projective
schemes}, Adv. Math. \textbf{109} (2), 228-287 (1994).

\bibitem[B]{Baez} J. Baez, \textit{The octonions}, Bull. Amer. Math. Soc.
\textbf{39}:2 (2002) 145-205; Errata, ibid. \textbf{42} (2005) 213; \texttt{%
math/0105155v4 [math.RA]}.

\bibitem[BFSS]{19} T. Banks, W. Fischler, S. H. Shenker and L. Susskind, $%
\mathit{M}$\textit{\ Theory as a Matrix Model: A Conjecture}, Phys. Rev.
\textbf{D55}, 5112-5128 (1997), \texttt{hep-th/9610043}.


\bibitem[BHQIT]{BHQIT} L. Borsten, M.J. Du , A. Marrani, W. Rubens, \textit{On the
Black-Hole/Qubit Correspondence}, Eur. Phys. J. Plus \textbf{126}, 37
(2011), \texttt{arXiv:1101.3559 [hep-th]}. L. Borsten, M.J. Duff, P. L\'{e}%
vay, \textit{The black-hole/qubit correspondence: an up-to-date review},
Class. Quant. Grav. \textbf{29}, 224008 (2012), \texttt{arXiv:1206.3166
[hep-th]}.

\bibitem[BKPPS]{beyond e11} G. Bossard, A. Kleinschmidt, J. Palmkvist, C. N.
Pope and E. Sezgin, \textit{Beyond E}$_{11}$, JHEP \textbf{1705} (2017) 020,
\texttt{arXiv:1703.01305 [hep-th]}.

\bibitem[Bou]{bour} N. Bourbaki \emph{Groups et Alg\`ebres de Lie}, Hermann, Paris 1968.

\bibitem[Bro]{brown} R. B. Brown, \textit{Groups of type E}$_{7}$, J. Reine
Angew. Math. \textbf{236}, 79 (1969).

\bibitem[BS]{12} B. Kim and A. Schwarz, \textit{Formulation of M(atrix)
model in terms of octonions} (unpublished).

\bibitem[Ca]{carter} Carter, R. W. \textit{Simple Groups of Lie Type},
Wiley-Interscience, New-York, 1989.

\bibitem[CDGL]{Montecarlo} G. Cossu, M. D'Elia, A. Di Giacomo, B. Lucini, C.
Pica, \textit{Confinement: G}$_{2}$\textit{\ group case}, PoSLAT2007, 296
(2007), \texttt{arXiv:0710.0481 [hep-lat]}.

\bibitem[CDS]{17} A. Connes, M. R. Douglas and A. Schwarz, \textit{%
Noncommutative geometry and matrix theory: compactification on tori}, JHEP
\textbf{9802}, 003 (1998), \texttt{hep-th/9711162}.

\bibitem[CJ]{CJ-1} E. Cremmer and B. Julia, \textit{The }$\mathcal{N}\mathit{%
=8}$\textit{\ Supergravity Theory. 1. The Lagrangian}, Phys. Lett. \textbf{%
B80}, 48 (1978). E. Cremmer and B. Julia, \textit{The }$\mathit{SO(8)}$%
\textit{\ Supergravity}, Nucl. Phys. \textbf{B159}, 141 (1979).

\bibitem[Co94]{15} A. Connes : \textit{\textquotedblleft Noncommutative
Geometry"}, Boston, Academic Press (1994).

\bibitem[Co95]{22} A. Connes, \textit{Noncommutative geometry and reality},
J. Math. Phys. \textbf{36}, 6194 (1995).

\bibitem[Co96]{23} A. Connes, \textit{Gravity coupled with matter and the
foundation of noncommutative geometry}, Comm. Math. Phys. \textbf{155}, 109
(1996).

\bibitem[CS]{f4} C. Chevalley and R. D. Schafer, \textit{The Exceptional
Simple Lie Algebras }$\mathit{F}_{4}$\textit{\ and }$\mathit{E}_{6}$, Proc.
Natl. Acad. Sci. U.S.A. \textbf{35} (2), 137-141 (1950).

\bibitem[Cy]{Octonions} A. Cayley, \textit{On Jacobi's elliptic functions,
in reply to the Rev. Brice Bronwin; and on quaternions}, Philosophical
Magazine \textbf{26}, 208-211 (1845).

\bibitem[Da]{spacelike-Refs} T. Damour and M. Henneaux, $\mathit{E}_{10}$%
\textit{, }$\mathit{BE}_{10}$\textit{\ and arithmetical chaos in superstring
cosmology}, Phys. Rev. Lett. \textbf{86}, 4749--4752 (2001), \texttt{%
hep-th/0012172}. T. Damour, M. Henneaux, B. Julia and H. Nicolai, \textit{%
Hyperbolic Kac-Moody algebras and chaos in Kaluza-Klein models}, Phys. Lett.
\textbf{B509}, 323--330 (2001), \texttt{hep-th/0103094}.

\bibitem[deGr]{graaf} W.A. de Graaf : \textquotedblleft \textit{Lie
Algebras: Theory and Algorithms"} (North-Holland Mathematical Library
\textbf{56}, Elsevier, Amsterdam, 2000).

\bibitem[DFLV]{Spinor} R. D'Auria, S. Ferrara, M.A. Lledo, V.S. Varadarajan,
\textit{Spinor algebras}, J.Geom.Phys. \textbf{40} (2001) 101-128, \texttt{hep-th/0010124}.

\bibitem[DG]{DG10} J. Distler and S. Garibaldi, \textit{There is no
\textquotedblleft Theory of Everything" inside }$\mathit{E}_{8}$, Comm.
Math. Phys. \textbf{298} (2010), 419, \texttt{arXiv:0905.2658 [math.RT]}.

\bibitem[EP2]{EP2} P. Truini, A. Marrani, M. Rios, \textit{Exceptional Periodicity and Magic Star Algebras. II : Gradings and {\em HT}-Algebras}; to appear.

\bibitem[EP3]{EP3} P. Truini, W. De Graaf, A. Marrani, \textit{Exceptional Periodicity and Magic Star Algebras. III : The Algebra $\fqn$ and the Derivations of {\em HT}-Algebras}; in preparation.

\bibitem[EP4]{EP4} P. Truini, W. De Graaf, A. Marrani, \textit{Exceptional Periodicity and Magic Star Algebras. IV : Cubic T-Algebras}; in preparation.

\bibitem[exp]{exp} R. Coldea, D.A. Tennant, E.M. Wheeler, E. Wawrzynska, D.
Prabhakaran, M. Telling, K. Habicht, P. Smibidl, and K. Kiefer, \textit{%
Quantum criticality in an Ising chain: experimental evidence for emergent E}$%
_{8}$\textit{\ symmetry}, Science \textbf{327}, 177 (2010). D. Borthwick and
S. Garibaldi, \textit{Did a 1-dimensional magnet detect a 248-dimensional
Lie algebra?}, Not. Amer. Math. Soc. \textbf{58}, 1055 (2011), \texttt{%
arXiv:1012.5407 [math-ph]}.

\bibitem[FGT]{FGT} A. Marrani, C.-X. Qiu, S.-Y. D. Shih, A. Tagliaferro, B.
Zumino, \textit{Freudenthal Gauge Theory}, JHEP \textbf{1303}, 132 (2013),
\texttt{arXiv:1208.0013 [hep-th]}.

\bibitem[FK]{minimal} S. Ferrara and R. Kallosh, \textit{Creation of Matter in
the Universe and Groups of Type E}$_{7}$, JHEP \textbf{1112}, 096 (2011),
\texttt{arXiv:1110.4048 [hep-th]}. S. Ferrara, R. Kallosh, and A. Marrani,
\textit{Degeneration of Groups of Type E}$_{7}$\textit{\ and Minimal
Coupling in Supergravity,} JHEP \textbf{1206}, 074 (2012), \texttt{%
arXiv:1202.1290 [hep-th]}.

\bibitem[GHMR]{Het} D. Gross, J. Harvey, E. Martinec and R. Rohm, \textit{%
Heterotic string}, Phys. Rev. Lett. \textbf{54} (6), 502-505 (1985).

\bibitem[GM61]{GM} M. Gell-Mann, \textit{The Eightfold Way: A Theory of
strong interaction symmetry}, Synchrotron Laboratory Report CTSL-20,
California Inst. of Tech., Pasadena (1961).

\bibitem[GM64]{quark-1} M. Gell-Mann, \textit{A Schematic Model of Baryons
and Mesons}, Phys. Lett. \textbf{8} (3), 214-215 (1964).

\bibitem[GRS]{unif} F. G\"{u}rsey, P. Ramond and P. Sikivie, \textit{A
universal gauge theory model based on }$\mathit{E}_{6}$, Phys. Lett. \textbf{%
B60} (2), 177-180 (1976).

\bibitem[GST]{MESGT} M. G\"{u}naydin, G. Sierra, P. K. Townsend, \textit{%
Exceptional Supergravity Theories and the Magic Square}, Phys. Lett. \textbf{%
B133} , 72 (1983). M. G\"{u}naydin, G. Sierra and P. K. Townsend, \textit{%
The Geometry of }$\mathcal{N}\mathit{=2}$\textit{\ Maxwell-Einstein
Supergravity and Jordan Algebras}, Nucl. Phys. \textbf{B242}, 244 (1984).


\bibitem[HKPW]{deco} K. Holland, P. Minkowski, M. Pepe and U. J. Wiese,
\textit{Exceptional confinement in G}$_{2}$\textit{\ gauge theory}, Nucl.
Phys. \textbf{B668}, 207 (2003), \texttt{hep-lat/0302023}.

\bibitem[HS]{14} G. T. Horowitz and L. Susskind, \textit{Bosonic }$\mathit{M}
$\textit{\ Theory}, J. Math. Phys. \textbf{42}, 3152 (2001), \texttt{%
hep-th/0012037}.

\bibitem[HT]{HT-1} C. Hull and P. K. Townsend, \textit{Unity of Superstring
Dualities}, Nucl. Phys. \textbf{B438}, 109 (1995), \texttt{hep-th/9410167}.

\bibitem[Hu]{hum} J.E. Humphreys: \textquotedblleft \textit{Introduction to
Lie Algebras and Representation Theory"} (Springer-Verlag New York Inc., New
York, 1972).

\bibitem[HW]{20} P. M. Ho and Y. S. Wu, \textit{Noncommutative Geometry and }%
$\mathit{D}$\textit{-branes}, Phys. Lett. \textbf{B398}, 52-60 (1997),
\texttt{hep-th/9611233}.

\bibitem[JWVN]{JWVN} P. Jordan, J. von Neumann and E. P. Wigner, \textit{On
an algebraic generalization of the quantum mechanical formalism}, Ann. Math.
\textbf{35} (1934) no. 1, 29--64.

\bibitem[Kac]{kac} V.G. Kac : \textquotedblleft \textit{Infinite Dimensional
Lie Algebras"}, third edition (Cambridge University Press, Cambridge, 1990).


\bibitem[KLR]{matrix} J.P. Keating, N.Linden and Z. Rudnick, \textit{Random
Matrix Theory, The exceptional Lie groups, and L-functions}, J. Phys.
\textbf{A36} no. 12, 2933 (special RMT volume) (2003).

\bibitem[KS]{D-branes} W. Krauth and M. Staudacher, \textit{Yang-Mills
integrals for orthogonal, symplectic and exceptional groups}, Nucl. Phys.
\textbf{B584}, 641 (2000), \texttt{hep-th/0004076}.

\bibitem[Li]{L} A. G. Lisi, \textit{An exceptionally simple theory of
everything}, \texttt{arXiv:0711.0770 [hep-th]}.

\bibitem[Lo]{loos1} O.Loos : \textit{\textquotedblleft Jordan Pairs"}, Lect.
Notes Math. \textbf{460}, (Springer, 1975).

\bibitem[Ma13]{Marrani-JP-5} S. Ferrara, A. Marrani, B. Zumino, \textit{%
Jordan Pairs, }$E_{6}$\textit{\ and }$U$\textit{-Duality in Five Dimensions}%
, J. Phys. \textbf{A46} (2013) 065402, \texttt{arXiv:1208.0347 [math-ph]}.

\bibitem[Ma14]{Marrani-Truini-1} A. Marrani and P. Truini, \textit{%
Exceptional Lie Algebras, }$\mathit{SU(3)}$\textit{\ and Jordan Pairs Part
2: Zorn-type Representations}, J. Phys. \textbf{A47} (2014) 265202, \texttt{%
arXiv:1403.5120 [math-ph]}.

\bibitem[Ma16]{Marrani-Truini-Interactions} A. Marrani and P. Truini,
\textit{Exceptional Lie Algebras at the very Foundations of Space and Time},
p Adic Ultra. Anal. Appl. \textbf{8} (2016) no.1, 68-86, \texttt{%
arXiv:1506.08576 [hep-th]}.

\bibitem[Ma17]{Marrani-Truini-2} A. Marrani and P. Truini, \textit{%
Sextonions, Zorn Matrices, and }$\mathfrak{e}_{7\frac{1}{2}}$, Lett. Math.
Phys. \textbf{107} (2017) no.10, 1859-1875, \texttt{arXiv:1506.04604
[math.RA]}.

\bibitem[Md]{18} J. Madore, \textit{Noncommutative Geometry for Pedestrians}%
, Lecture given at the International School of Gravitation, Erice, \texttt{%
gr-qc/9906059}.


\bibitem[MPRR]{Magic-Non-Susy} A. Marrani, G. Pradisi, F. Riccioni and L.
Romano, \textit{Non-Supersymmetric Magic Theories and Ehlers Truncations},
Int. J. Mod. Phys. \textbf{A32} (2017) no.19-20, 1750120, \texttt{%
arXiv:1701.03031 [hep-th]}.

\bibitem[MR]{Romano} A. Marrani and L. Romano, \textit{Orbits in
Non-Supersymmetric Magic Theories}, Int. J. Mod. Phys. A 2019 (in press),
\texttt{arXiv:1906.05830 [hep-th]}.

\bibitem[MS]{Magic Square} H. Freudenthal, \textit{Beziehungen der }$\mathit{%
E}_{7}$\textit{\ und }$\mathit{E}_{8}$\textit{\ zur oktavenebene I-II},
Nederl. Akad. Wetensch. Proc. Ser. \textbf{57} (1954) 218--230. J. Tits,
\textit{Interpr\'{e}tation g\'{e}ometriques de groupes de Lie simples
compacts de la classe E}, M\'{e}m. Acad. Roy. Belg. Sci \textbf{29} (1955)
3. H. Freudenthal, \textit{Beziehungen der E}$_{7}$\textit{\ und E}$_{8}$%
\textit{\ zur oktavenebene IX}, Nederl. Akad. Wetensch. Proc. Ser. \textbf{%
A62} (1959) 466--474. B. A. Rosenfeld, \textit{Geometrical interpretation of
the compact simple Lie groups of the class E}, (Russian) Dokl. Akad. Nauk.
SSSR \textbf{106} (1956) 600--603. J. Tits, \textit{Alg\`{e}bres
alternatives, alg\`{e}bres de Jordan et alg\`{e}bres de Lie exceptionnelles}%
, Indag. Math. \textbf{28} (1966) 223--237.

\bibitem[Mu]{Mukai} S. Mukai, \textit{Simple Lie algebra and Legendre variety%
}, Nagoya S\={u}ri Forum, \textbf{3} (1996), 1-12.

\bibitem[Ni]{extended-Refs} H. Nicolai, \textit{The integrability of }$%
\mathcal{N}\mathit{=16}$\textit{\ supergravity}, Phys. Lett. \textbf{B194},
402 (1987). R. P. Geroch, \textit{A method for generating solutions of
Einstein's equations}, J. Math. Phys. \textbf{12}, 918--924 (1971). R. P.
Geroch, \textit{A method for generating new solutions of Einstein's
equation. 2}, J. Math. Phys. \textbf{13}, 394--404 (1972). B. Julia, in :
\textit{\textquotedblleft Lectures in applied mathematics"}, vol. \textbf{21}%
, p. 355. AMS-SIAM, 1985. B. Julia, \textit{Group Disintegrations}. Invited
paper presented at Nuffield Gravity Workshop, Cambridge, Eng., Jun 22 - Jul
12, 1980. S. Mizoguchi, \textit{E}$_{10}$\textit{\ symmetry in
one-dimensional supergravity}, Nucl. Phys. \textbf{B528}, 238--264 (1998),
\texttt{hep-th/9703160}.

\bibitem[PR]{10} T. Pengpan and P. Ramond, \textit{M(ysterious) Patterns in }%
$\mathit{SO(9)}$, Phys. Rept. \textbf{315}, 137-152 (1998), \texttt{%
hep-th/9808190}.

\bibitem[Ra]{11} P. Ramond, \textit{Exceptional Groups and Physics}, Plenary
Talk delivered at the Conference Groupe 24, Paris, July 2002, \texttt{%
arXiv:hep-th/0301050v1}.

\bibitem[RMC19]{Chester-1} M. Rios, A. Marrani, D. Chester, \textit{Geometry
of exceptional super Yang-Mills theories}, Phys. Rev. \textbf{D99} (2019)
no.4, 046004, \texttt{arXiv:1811.06101 [hep-th]}.

\bibitem[RMC19-2]{Chester-2} M. Rios, A. Marrani, D. Chester, \textit{%
Exceptional Super Yang-Mills in }$\mathit{27+3}$\textit{\ and Worldvolume
M-Theory}, \texttt{arXiv:1906.10709 [hep-th]}.

\bibitem[Sa09]{Sati-1} H. Sati, $\mathbb{OP}^{2}$ \textit{Bundles in }$%
\mathit{M}$\textit{-Theory}, Commun. Num. Theor. Phys. \textbf{3} (2009)
495, \texttt{arXiv:0807.4899 [hep-th]}.

\bibitem[Sa11]{Sati-2} H. Sati, \textit{On the geometry of the
supermultiplet in }$\mathit{M}$\textit{-theory}, Int. J. Geom. Meth. Mod.
Phys. \textbf{8} (2011) 1, \texttt{arXiv:0909.4737 [hep-th]}.


\bibitem[Sm]{13} L. Smolin, \textit{The exceptional Jordan algebra and the
matrix string}, \texttt{hep-th/0104050}.


\bibitem[Sz]{21-bis} R. J. Szabo, $\mathit{D}$\textit{-Branes, Tachyons and }%
$\mathit{K}$\textit{-Homology}, Mod. Phys. Lett. \textbf{A17}, 2297 (2002),
\texttt{hep-th/0209210}.


\bibitem[Tr11]{Truini} P. Truini, \textit{Exceptional Lie Algebras, }$\mathit{%
SU(3)}$\textit{\ and Jordan Pairs}, Pacific J. Math. \textbf{260}, 227
(2012), \texttt{arXiv:1112.1258 [math-ph]}.

\bibitem[Tr19]{Truini-EUG} P. Truini, \textit{Vertex operators for an
expanding universe}, invited paper in {\it Symmetries and Order: Algebraic Methods in Many Body Systems}, Yale 5-6 October 2018, in honor of Francesco Iachello, on the occasion of his retirement; \texttt{arXiv:1901.07916 [physics.gen-ph]}.

\bibitem[Tr19-2]{Truini-QGR-forthcoming} P. Truini \textit{et al.}; to appear.

\bibitem[TRM17]{Mile-High-EP} P. Truini, M. Rios, A. Marrani, \textit{The
Magic Star of Exceptional Periodicity},  Contemp. Math. \textbf{721} (2019),
277-297, \texttt{arXiv:1711.07881 [hep-th]}.

\bibitem[TRM18]{Group32-1} A. Marrani, P. Truini, M. Rios, \textit{The Magic
of Being Exceptional}, J. Phys. Conf. Ser. \textbf{1194} (2019) no.1,
012075, \texttt{arXiv:1811.11208 [hep-th]}.

\bibitem[TRM18-2]{Group32-2} P. Truini, A. Marrani, M. Rios, \textit{Magic
Star and Exceptional Periodicity: an approach to Quantum Gravity}, J. Phys.
Conf. Ser. \textbf{1194} (2019) no.1, 012106, \texttt{arXiv:1811.11202
[hep-th]}.

\bibitem[TW]{West} A. G. Tumanov and P. West, $\mathit{E}_{11}$\textit{\ in
11D}, Phys. Lett. \textbf{B758} (2016) 278, \texttt{1601.03974 [hep-th]}.

\bibitem[Vi]{Vinberg} E.B. Vinberg, \textit{The theory of Convex Homogeneous
Cones}, in Transaction of the Moscow Mathematical Society for the year 1963,
340-403, American Mathematical Society, Providence RI 1965.

\bibitem[Vo]{KLV} D. Vogan, \textit{The character table for E}$_{8}$,
Notices of the AMS \textbf{54} (2007), no. 9, 1022.

\bibitem[Wi]{9} E. Witten, \textit{String theory dynamics in various
dimensions}, Nuclear Physics \textbf{B443} (1), 85-126 (1995), \texttt{%
hep-th/9503124}.

\bibitem[Zw64]{quark-2} G. Zweig, \textit{An }$\mathit{SU(3)}$\textit{\
Model for Strong Interaction Symmetry and its Breaking: II}, CERN Report No.
8419/TH.401 (1964).

\end{thebibliography}
\end{document}